\documentclass{article}
\usepackage{graphicx,bm,url,amssymb,amsmath,hyperref,color}
\usepackage{geometry}

\def\no{\noindent}
\def\pmatrix{\left(\begin{array}}
\def\endpmatrix{\end{array}\right)}

\def\blue{}

\def\NN{\mathbb{N}}

\def\RR{\mathbb{R}}

\def\I{{\cal I}}

\def\P{{\cal P}}

\def\h{{\bm{h}}}

\def\dd{\mathrm{d}}

\def\diag{\mathrm{diag}}

\newtheorem{theo}{Theorem}

\newtheorem{cor}{Corollary}
\newtheorem{rem}{Remark}

\def\proof{\noindent\underline{Proof}\quad}
\def\QED{\mbox{$\hfill{}~ \Box$}}

\def\bfc{{\bm{c}}}

\def\bff{{\bm{f}}}

\def\bfh{{\bm{h}}}

\def\bfzero{{\bm{0}}}
\def\bfuno{{\bm{1}}}

\def\bfalfa{{\bm{\alpha}}}

\def\bfgamma{{\bm{\gamma}}}
\def\bfdelta{{\bm{\delta}}}

\def\bfomega{{\bm{\omega}}}

\def\aa{{\alpha}}

\def\cpr{$^\copyright\,$}

\def\oP{{\pi}}
\def\FDEim{{\tt fde\_pi2\_im}}
\def\FDEpc{{\tt fde\_pi12\_pc}}
\def\FDEpcc{{\tt fde\_pi12\_pc-10}}

\def\taa{{\tilde\aa}}

\def\JP{Jacobi-Pi\~neiro\,}

\def\fhbvmnew{{\tt fhbvm2\underline{~}2}}

\begin{document}

\title{\blue A Multi-Order Extension of Fractional HBVMs (FHBVMs)}

\author{
Luigi Brugnano\,\footnote{Dipartimento di Matematica e Informatica ``U.\,Dini'', 
             Universit\`a di Firenze,  Italy,~    \url{luigi.brugnano@unifi.it}} 
   \and  
Gianmarco Gurioli\,\footnote{Dipartimento di Matematica e Informatica ``U.\,Dini'',
             Universit\`a di Firenze,  Italy, ~ \url{gianmarco.gurioli@unifi.it}} 
     \and 
Felice Iavernaro\,\footnote{Dipartimento di Matematica, 
           Universit\`a di Bari ``Aldo Moro'', Italy, ~ \url{felice.iavernaro@uniba.it}} 
     \and  
Mikk Vikerpuur\,\footnote{Institute of Mathematics and Statistics, 
           University of Tartu, Estonia, ~ \url{mikk.vikerpuur@ut.ee}} 
            }

\maketitle

\begin{abstract} 
The efficient numerical solution of fractional differential equations has been recently tackled through the definition of Fractional HBVMs (FHBVMs), a class of Runge-Kutta type methods.  Corresponding Matlab\cpr codes have been also made available on the internet, proving to be very competitive w.r.t. existing ones. However, so far, FHBVMs have been given for solving systems of fractional differential equations with the same order of fractional derivative, whereas the numerical solution of multi-order problems (i.e., problems in which different orders of fractional derivatives occur) has not been handled, yet. Due to their relevance in applications, in this paper we propose an extension of FHBVMs for addressing fractional multi-order  problems,  providing full details for such an approach. A corresponding Matlab\cpr code, handling the case of two different fractional orders, is also made available, proving very effective for numerically solving these problems.

\medskip
\no{\bf Keywords:} Fractional Differential Equations, FDEs, Caputo derivative, Fractional HBVMs, FHBVMs, fractional multi-order  problems.

\medskip
\no{\bf MSC:}  34A08, 65R20, 65-04.
\end{abstract}

\section{Introduction} 
The idea of extending the definition of derivatives to the case of a non-integer order of differentiation dates back to 1695, to the correspondence between Leibniz and de l’H\^{o}pital \cite{Pod99}, when the meaning of the derivative of order one half was discussed. Although its theoretical foundation is centuries old, the concept of fractional derivative has recently evolved from a purely mathematical curiosity into a vital tool for modeling across a variety of scientific fields (see, e.g., \cite{Hen08,Koh90,Magi10,NigArb07,Petras2011,SZBCC2018,Ucha22}, to mention a few), due to its ability to encompass memory terms. In fact, unless the ODE case, the differential operator here considered is of a nonlocal nature \cite{Diethelm2010,Pod99}. For this reason, the numerical solution of FDEs has been the subject of many researches in the last years, starting from the pioneering work of Lubich \cite{Lub85}, where fractional versions of linear multistep formulae are considered: we refer, e.g., to \cite{Die04,Die05,FPV2023,Garr11,Garr18,Li16,PTV2017}, to mention a few of them. More recently, the class of {\em Fractional HBVMs (FHBVMs)} \cite{BBBI2024} has been introduced for the efficient numerical solution of initial value problems of fractional differential equations (FDE-IVPs) in the Caputo form, with one fractional derivative. 

{\blue
\begin{rem}\label{hbvms} 
The acronym HBVMs stands for {\em ``Hamiltonian Boundary Value Methods''}, a class of energy-conserving Runge-Kutta methods for the efficient numerical solution of Hamiltonian ODE problems (see, e.g., the monograph  \cite{book}). As matter of fact, when the fractional order of the equation becomes 1 (i.e., the FDE reduces to an ODE), FHBVMs do coincide with HBVMs.
\end{rem}
}

The basic steps leading to the definition of FHBVMs are the following ones:
\begin{itemize}
\item local expansion of the vector field along  a suitable Jacobi polynomial basis, naturally induced by the {\blue fractional problem at hand (when the equation is an ordinary one, the Jacobi basis reduces to the Legendre polynomial basis, thus resulting into HBVMs)};

\item discretization of the corresponding Fourier coefficients, by using a corresponding Gauss-Jacobi quadrature;

\item solution of the generated discrete problems, defined on a  suitable discrete mesh.

\end{itemize} 
FHBVMs have been efficiently implemented in the two Matlab\cpr codes:
\begin{itemize}
\item {\tt fhbvm}, \cite{BGI2024,BGI2025} using either a graded or a uniform mesh (depending on the features of the problem at hand);
\item {\tt fhbvm2}, \cite{BGIV2025} possibly also using a mixed mesh (initially graded and subsequently uniform).
\end{itemize}
An extensive experimentation \cite{BGIV2025_1} has proved that the above codes are extremely competitive w.r.t. publicly available ones, based on different numerical methods, due to the fact that the implemented FHBVM methods can achieve, alike in the ODE case \cite{ABI2020}, a spectrally accurate solution.

Nevertheless, multi-order systems of equations, i.e., systems where different fractional orders appear, are of interest in a number of applications, including material science \cite{Torvik1984,Pod99}), population dynamics \cite{Arciga2019}, epidemic modelling \cite{Chen2021,Li2019}, pharmacokinetics \cite{Vero2010}, and chaotic dynamics \cite{Baleanu2021,Grigorenko2003,Li2015,Sheu2006}. {\blue
In this respect, the extension of FHBVMs for such problems poses important challenges.  To mention a few of them: the handling of different Jacobi quadratures; a more involved structure and solution of the generated discrete problems; a more involved theoretical analysis of the methods. Their study is precisely the goal of this paper.}

Consequently, we shall hereafter consider multi-order problems in the following general form, for \,$\aa_i\in(\ell-1,\ell)$,\, $i=1\dots,\nu$,\, $\ell\in\NN$, $\ell\ge1$:
\begin{eqnarray}\nonumber
y_i^{(\aa_i)}(t) &=& f_i(y(t)), \qquad t\in[0,T], \\[1mm] \nonumber
y_i^{(\iota)}(0) &=& y_{i0}^\iota\in\RR^{m_i},\qquad \iota=0,\dots,\ell-1, \qquad i = 1,\dots,\nu, \\
y &=& (y_1^\top,\dots,y_\nu^\top)^\top \in\RR^m, \qquad m = \sum_{i=1}^\nu m_i.\label{fde1}
\end{eqnarray}
Here, for all \,$i=1,\dots,\nu$, ~assuming ~$y_i\in A^\ell([0,T])$ ~(i.e., $y_i^{(\ell-1)}$ is absolutely continuous),
\begin{equation}\label{caputo}
y_i^{(\aa_i)}(t) ~:=~ \frac{1}{\Gamma(\ell-\aa_i)}\int_0^t (t-x)^{\ell-\aa_i-1}y_i^{(\ell)}(x)\dd x
\end{equation}
is the Caputo derivative of $y_i$ and, for the sake of brevity, in (\ref{fde1}) we have omitted $t$ as a formal argument of $f_i$. 
The approach described in the following sections reduces to the usual one defining FHBVMs, when the $\{\aa_i\}$ are all equal (or, equivalently, $\nu=1$), but differs significantly both in the metodology and the implementations strategies, when they are different from each other. 

{\blue While the numerical solution of linear or nonlinear multi-term FDEs\,\footnote{\blue I.e., FDEs with different fractional derivatives appearing in the equation.} involving Caputo-type fractional derivatives has, over the years, attracted considerable attention by a large number of authors, the general multi-order problem in the form (\ref{fde1}), has received significantly less study. Nevertheless, we can point the reader interested in the numerical analysis of such problems to the Adomian decomposition and the variational iteration method studied in \cite{Momani2007}, the differential transform technique used in \cite{Erturk2008}, the Chebyshev matrix approach implemented in \cite{Atabakzadeh2013, Dabiri2017}, the application of the fractional Jacobi-Picard iteration method in \cite{Ansari2024}, wavelet-type solutions in \cite{Chen2015,Wang2019}, and the trapezoidal and product integration methods in \cite{Garr18}. We also note that there exist some prior numerical methods implementing a spectral-type approach, with basis functions consisting of fractional-order Jacobi functions \cite{Bhrawy2016,Faghih2022} or Chelyshkov polynomials \cite{Ahmed2022}.
}

With this premise, the structure of the paper is as follows. At first, in Section~\ref{fhbvm_sec} we explain the possible approaches for extending FHBVMs to handle problem (\ref{fde1}). Later on, in Sections~\ref{nugt1} and \ref{mops}, we provide full details for the corresponding discretization procedure, emphasizing the major difference in the choice of the discrete abscissae and quadratures. The practical solution of the generated discrete problems is then detailed in Section~\ref{newton}, {\blue and the stability and accuracy analysis of the methods is given in Section~\ref{analmet}. }A new version of the code {\tt fhbvm2} in \cite{BGIV2025}, named \fhbvmnew, is {\blue then presented} in Section~\ref{code}, handling at the moment the case $\nu=2$. This latter code is later used in the numerical tests  reported in Section~\ref{numtest}. At last, some conclusions are given in Section~\ref{fine}. 

\section{Fractional HBVMs (FHBVMs)}\label{fhbvm_sec}
To begin with, we recall that, under suitable regularity assumptions on the vector fields $f_i$, the solution of problem (\ref{fde1}) can be formally written  (see, e.g., \cite{DGS2020}) as:
\begin{eqnarray}\nonumber
y_i(t) &=& \sum_{\iota=0}^{\ell-1} \frac{y_{i0}^\iota}{\iota!}t^\iota  ~+~ \frac{1}{\Gamma(\aa_i)}\int_0^t (t-x)^{\aa_i-1} y_i^{(\aa_i)}(x)\dd x\\
\nonumber
         &=& \sum_{\iota=0}^{\ell-1} \frac{y_{i0}^\iota}{\iota!}t^\iota  ~+~ \frac{1}{\Gamma(\aa_i)}\int_0^t (t-x)^{\aa_i-1} f_i(y(x))\dd x\\[2mm]
         &\equiv& T_{i\ell}(t) + I^{\aa_i} f_i(y(t)), \qquad t\in[0,T], \qquad i=1,\dots,\nu,\label{sol}
\end{eqnarray}
with $y$ defined according to (\ref{fde1}), and
$$I^{\aa_i}f_i(y(t)) = \frac{1}{\Gamma(\aa_i)}\int_0^t (t-x)^{\aa_i-1}f(_i(y(x))\dd x$$
the Riemann-Liouville integral corresponding to (\ref{caputo}).

Let us now consider, for explanation purposes, the case where problem (\ref{fde1}) is, at first, solved over the reference interval $[0,h]$. To this end, we use the families of Jacobi polynomials $\{P_j^i\}_{j\ge0}$, $i=1,\dots,\nu$, respectively orthonormal w.r.t. the weighting functions:
 \begin{equation}\label{ortoi}
\omega_i(c) = \aa_i(1-c)^{\aa_i-1}, \quad c\in[0,1], \quad \Longrightarrow \quad \int_0^1 \omega_i(c) P_r^i(c)P_j^i(c)\dd c = \delta_{rj},\quad r,j=0,1,\dots,
\end{equation}
with, as is usual, $\delta_{rj}$ denoting the Kronecker symbol. Next, we expand each vector field $f_i$ in (\ref{fde1}) along the corresponding family of polynomials, thus obtaining the equivalent problem
\begin{eqnarray}\nonumber
y_i^{(\aa_i)}(ch) &=& f_i(y(ch)) ~\equiv~ \sum_{j\ge0} P_j^i(c)\gamma_{ij}(y), \qquad c\in[0,1],\\ \label{fde1i}
y_i^{(\iota)}(0) &=& y_{i0}^\iota\in\RR^{m_i},\qquad \iota=0,\dots,\ell-1,\qquad i=1,\dots,\nu,
\end{eqnarray}
with the corresponding Fourier coefficients defined as
\begin{equation}\label{gamji}
\gamma_{ij}(y) = \int_0^1 \omega_i(\tau)P_j^i(\tau)f_i(y(\tau h))\dd\tau, \qquad j=0,1,\dots.
\end{equation}
Integrating (\ref{fde1i}) side by side, and imposing the initial conditions, one then obtains
\begin{equation}\label{sol2}
y_i(ch) = T_{i\ell}(ch) + h^{\aa_i}\sum_{j\ge0} I^{\aa_i} P_j^i(c)\gamma_{ij}(y), \qquad c\in[0,1],\qquad i=1,\dots,\nu,
\end{equation}
which is equivalent to (\ref{sol}), for $t\in[0,h]$. Moreover, by considering that
$$P_0^i(c)\equiv 1 \qquad \mbox{and}\qquad I^{\aa_i} P_j^i(1)=\frac{\delta_{j0}}{\Gamma(\aa_i+1)},$$ 
from (\ref{ortoi}), (\ref{gamji}), and (\ref{sol2}) one derives:
\begin{eqnarray}\nonumber
y_i(h) &=& T_{i\ell}(h) + \frac{h^{\aa_i}}{\Gamma(\aa_i+1)}\gamma_{i0}(y) ~=~T_{i\ell}(h) + \frac{h^{\aa_i}}{\Gamma(\aa_i)\aa_i}\gamma_{i0}(y)
\\[1mm] \nonumber
&=& T_{i\ell}(h) + \frac{1}{\Gamma(\aa_i)}\int_0^h (h-x)^{\aa_i-1}f_i(y(x))\dd x 
~\equiv~  T_{i\ell}(h) + I^{\aa_i} f_i(y(h)),\label{solh}
\end{eqnarray}
that is, (\ref{sol}) at $t=h$. 
Polynomial approximations of degree $s-1$ to (\ref{fde1i}) are obtained by truncating the infinite series  to finite sums with $s$ terms:
\begin{eqnarray}\nonumber
\sigma_i^{(\aa_i)}(ch) &=& \sum_{j=0}^{s-1} P_j^i(c)\gamma_{ij}(\sigma), \qquad c\in[0,1],\\ \nonumber 
\sigma_i^{(\iota)}(0) &=& y_{i0}^\iota\in\RR^{m_i},\qquad \iota=0,\dots,\ell-1, \qquad i=1,\dots,\nu,\\[2mm]
\sigma &=& (\sigma_1^\top,\dots,\sigma_\nu^\top)^\top\in\RR^m,\label{sigalfa}
\end{eqnarray}
with the Fourier coefficients $\gamma_{ij}(\sigma)$ defined according to (\ref{gamji}), by formally replacing $y$ with $\sigma$. As is clear, this is equivalent to requiring the residuals be orthogonal to all polynomials of degree $s-1$, w.r.t. the inner product defined by (\ref{ortoi}):\footnote{This feature is common to HBVMs, obtained when $\aa_i=1$ \cite{ABI2023,book}.}
\begin{equation}\label{spectacc}
\int_0^1\omega_i(c)P_j^i(c)\left[ \sigma_i^{(\aa_i)}(ch)-f_i(\sigma_i(ch))\right]\dd c = 0, \qquad j=0,\dots,s-1, \qquad i=1,\dots,\nu.
\end{equation}
Similarly to (\ref{sol2}), the approximations $\sigma_i$ are then given by integrating (\ref{sigalfa}) side by side and imposing the initial conditions:
\begin{eqnarray}\nonumber
\sigma_i(ch) &=&T_{i\ell}(ch) + \frac{h^{\aa_i}}{\Gamma(\aa_i)}\int_0^c (c-\tau)^{\aa_i-1}\sigma_i^{(\aa_i)}(\tau h)\dd\tau
~\equiv~ T_{i\ell}(ch)+I^{\aa_i} \sigma_i^{(\aa_i)}(ch)\\
&=& T_{i\ell}(ch) + h^{\aa_i}\sum_{j=0}^{s-1} I^{\aa_i} P_j^i(c)\gamma_{ij}(\sigma), \qquad c\in[0,1], \qquad i=1,\dots,\nu.
\end{eqnarray}
In particular, at $t=h$ one derives, similarly as in (\ref{solh}):
\begin{eqnarray*} 
\sigma_i(h) &=& T_{i\ell}(h) +\frac{h^{\aa_i}}{\Gamma(\aa_i+1)}\gamma_{i0}(\sigma) 
~=~ T_{i\ell}(h) + \frac{1}{\Gamma(\aa_i)}\int_0^h (h-x)^{\aa_i-1}f_i(\sigma(x))\dd x \\[2mm]
&\equiv& T_{i\ell}(h) + I^{\aa_i} f_i(\sigma(h)), \qquad i=1,\dots,\nu.
\end{eqnarray*} 

Next, let us return to the interval $[0,T]$, and consider the case where a mesh
\begin{equation}\label{hn}
t_n = t_{n-1}+h_n, \qquad n=1,\dots,N, \qquad t_0 = 0, \qquad  t_N=T,
\end{equation}
is given. Consequently, now $\sigma(t)$ becomes a piecewise approximation to $y(t)$ on the whole interval.
Moreover, let us denote:
\begin{eqnarray}\nonumber
\sigma_{in}(ch_n) &=& \sigma_i(t_{n-1}+ch_n), \\[1mm] \nonumber 
\sigma_{in}^{(\aa_i)}(ch_n) &=& \sum_{j=0}^{s-1} P_j^i(c) \gamma_{ij}(\sigma^n),\\[1mm]\nonumber
\gamma_{ij}(\sigma^n) &=& \int_0^1 \omega_i(c)P_j^i(c)f_i(\sigma^n(ch_n))\dd c,\qquad j=0,\dots,s-1,\\[1mm]\nonumber
T_{i\ell}^n(ch_n) &=& T_{i\ell}(t_{n-1}+ch_n), \qquad i=1,\dots,\nu, \\[2mm]
\sigma^n(ch_n) &=&\sigma(t_{n-1}+ch_n), \qquad c\in[0,1], \qquad n=1,\dots,N.\label{sigin}
\end{eqnarray}
Consequently, for $t=t_{n-1}+ch_n\in[t_{n-1},t_n]$, $c\in[0,1]$, one has:
\begin{eqnarray}\nonumber
\sigma_i(t) &\equiv& \sigma_{in}(ch_n)~=~ T_{i\ell}(t) \,+\, \frac{1}{\Gamma(\aa_i)}\int_0^t (t-x)^{\aa_i-1}\sigma_i^{(\aa_i)}(x)\dd x\\
\nonumber
&=& T_{i\ell}(t_{n-1}+ch_n) + \frac{1}{\Gamma(\aa_i)}\int_0^{t_{n-1}+ch_n} (t_{n-1}+ch_n-x)^{\aa_i-1}\sigma_i^{(\aa_i)}(x)\dd x \\
\nonumber
&=& T_{i\ell}^n(ch_n) \,+\, \frac{1}{\Gamma(\aa_i)}\left[ \sum_{\mu=1}^{n-1} \int_{t_{\mu-1}}^{t_\mu} (t_{n-1}+ch_n-x)^{\aa_i-1}\sigma_i^{(\aa_i)}(x)\dd x \right.\\ 
\nonumber
&& \qquad\left. +\int_{t_{n-1}}^{t_{n-1}+ch_n}(t_{n-1}+ch_n-x)^{\aa_i-1}\sigma_i^{(\aa_i)}(x)\dd x\right]\\
\nonumber
&=& T_{i\ell}^n(ch_n) \,+\, \left[ \sum_{\mu=1}^{n-1} \frac{h_\mu^{\aa_i}}{\Gamma(\aa_i)}\int_0^1 \left(\frac{t_{n-1}+ch_n-t_{\mu-1}}{h_\mu}-\tau\right)^{\aa_i-1}\sigma_{i\mu}^{(\aa_i)}(\tau h_\mu)\dd \tau \right.\\ 
\nonumber
&& \qquad\left. +\,\frac{h_n^{\aa_i}}{\Gamma(\aa_i)}\int_0^c(c-\tau)^{\aa_i-1}\sigma_{in}^{(\aa_i)}(\tau h_n)\dd \tau\right]
\\ \nonumber
&=& T_{i\ell}^n(ch_n) \,+\, \left[ \sum_{\mu=1}^{n-1} \frac{h_\mu^{\aa_i}}{\Gamma(\aa_i)}\int_0^1 \left(\frac{t_{n-1}+ch_n-t_{\mu-1}}{h_\mu}-\tau\right)^{\aa_i-1}\,\sum_{j=0}^{s-1} P_j^i(\tau) \gamma_{ij}(\sigma^\mu)\dd \tau \right.\\ 
\nonumber
&& \qquad\left. +\,
\frac{h_n^{\aa_i}}{\Gamma(\aa_i)}\int_0^c(c-\tau)^{\aa_i-1}\sum_{j=0}^{s-1} P_j^i(\tau) \gamma_{ij}(\sigma^n)\dd \tau\right]\\
\nonumber
&=& T_{i\ell}^n(ch_n) \,+\,  \sum_{\mu=1}^{n-1} h_\mu^{\aa_i}\sum_{j=0}^{s-1} J_j^i\left(\frac{t_{n-1}+ch_n-t_{\mu-1}}{h_\mu}\right)\gamma_{ij}(\sigma^\mu)  \,+\, h_n^{\aa_i}\sum_{j=0}^{s-1} I^{\aa_i}P_j^i(c) \gamma_{ij}(\sigma^n)\\
&\equiv& \phi_i^n(ch_n) \,+\, h_n^{\aa_i}\sum_{j=0}^{s-1} I^{\aa_i}P_j^i(c) \gamma_{ij}(\sigma^n),\qquad c\in[0,1], \qquad i=1,\dots,\nu,
\label{sigin1}
\end{eqnarray}
having set, for $x\ge1$,
\begin{equation}\label{Jji}
J_j^i(x) := \frac{1}{\Gamma(\aa_i)}\int_0^1(x-\tau)^{\aa_i-1}P_j^i(\tau)\dd\tau, \qquad j=0,\dots,s-1, \qquad i=1,\dots,\nu,
\end{equation}
and the {\em memory} term
\begin{equation}\label{fijin}
\phi_i^n(ch_n) := T_{i\ell}^n(ch_n) \,+\, \sum_{\mu=1}^{n-1} h_\mu^{\aa_i}\sum_{j=0}^{s-1} J_j^i\left(\frac{t_{n-1}+ch_n-t_{\mu-1}}{h_\mu}\right)\gamma_{ij}(\sigma^\mu), \quad c\in[0,1], \qquad i=1,\dots,\nu.
\end{equation}
\begin{rem} For the efficient numerical evaluation of the integrals $I^{\aa_i}P_j^i(c)$ in (\ref{sol2}), and $J_j^i(x)$ in (\ref{Jji}), we refer to \cite[Section\,3.2]{BGI2024}.\end{rem}

It must be emphasized that, since the memory term only depends on the past integration steps, according to (\ref{sigin})--(\ref{fijin}) the discrete problem at step $n$ consists in finding the $s\nu$ Fourier coefficients:
\begin{equation}
\gamma_{ij}(\sigma^n) = \int_0^1 \omega_i(c)P_j^i(c)f_i(\sigma^n(ch_n))\dd c,\qquad j=0,\dots,s-1,\qquad i=1,\dots,\nu.\label{stepn0}
\end{equation}
In more details, by setting hereafter $I_r\in\RR^{r\times r}$ the identity matrix, $$\bfuno_r=\left(1,\,\dots,\,1\right)^\top\in\RR^r,
\qquad \bfalfa = (\aa_1,\dots,\aa_\nu),$$
introducing the vectors, with $\otimes$ denoting the Kronecker product, 
\begin{eqnarray}\nonumber
\gamma_i(\sigma^n) &=& \pmatrix{c} \gamma_{i0}(\sigma^n)\\ \vdots\\ \gamma_{i,s-1}(\sigma^n)\endpmatrix\in\RR^{s m_i}, \qquad i=1,\dots,\nu,\\[2mm]
 \label{notations0}
\bfgamma(\sigma^n) &=& \pmatrix{c} \gamma_1(\sigma^n)\\ \vdots\\ \gamma_\nu(\sigma^n)\endpmatrix\in\RR^{sm},
\end{eqnarray}
the matrices and functions
\begin{eqnarray}\nonumber
P^i(c) &=& \pmatrix{ccc} P_0^i(c)\\ &\ddots\\ &&P_{s-1}^i(c)\endpmatrix\otimes I_{m_i}\in\RR^{sm_i \times s m_i}, \qquad i=1,\dots,\nu,\\[2mm]
\nonumber
P(c)   &=& \pmatrix{ccc} P^1(c)\\ & \ddots\\ && P^\nu(c)\endpmatrix\in\RR^{sm\times sm},
\\[2mm] \nonumber 
I^{\aa_i}P^i(c) &=&  \pmatrix{ccc} I^{\aa_i}P_0^i(c), &\dots,& I^{\aa_i}P_{s-1}^i(c)\endpmatrix\otimes I_{m_i}\in\RR^{m_i\times s m_i}, \qquad i=1,\dots,\nu, \\[2mm]
\nonumber
h_n^{\bfalfa} &=& \pmatrix{ccc} h_n^{\aa_1}I_{m_1}\\ &\ddots\\ &&h_n^{\aa_\nu}I_{m_\nu}\endpmatrix\in\RR^{m\times m},\\[2mm]
\nonumber
I^{\bfalfa}P(c)   &=& \pmatrix{ccc} I^{\aa_1}P^1(c)\\ &\ddots\\ && I^{\aa_\nu}P^\nu(c)\endpmatrix\in\RR^{m\times sm},\\[2mm]
\label{notations1}
\phi^n(ch_n) &=&\pmatrix{c} \phi_1^n(ch_n)\\ \vdots \\ \phi_\nu^n(ch_n)\endpmatrix\in\RR^m,
\end{eqnarray}
and taking into account that 
\begin{equation}\label{sigman}
\sigma^n(ch_n) = \phi^n(ch_n) + h_n^{\bfalfa}I^{\bfalfa} P(c)\bfgamma(\sigma^n),
\end{equation}
one has that (\ref{stepn0}) can be rewritten as:
\begin{eqnarray}\nonumber
\gamma_{ij}(\sigma^n) &=& \int_0^1 \omega_i(c)P_j^i(c)f_i\left(\phi^n(ch_n)+h_n^\bfalfa I^\bfalfa P(c)\bfgamma(\sigma^n)\right)\dd c,\\[1mm]
&& \qquad j=0,\dots,s-1,\qquad i=1,\dots,\nu. \label{stepn}
\end{eqnarray}
This problem, in turn, can be cast in vector form as
\begin{equation}\label{stepn_v}
\bfgamma(\sigma^n) = \int_0^1 \bfomega(c) P(c)  \bff\left( \phi^n(ch_n) + h_n^{\bfalfa}I^{\bfalfa} P(c)\bfgamma(\sigma^n)\right)\dd c,
\end{equation}
having set
\begin{eqnarray}\nonumber
\bff(\sigma) &=& \pmatrix{c} f_1(\sigma)\\ \vdots\\ f_\nu(\sigma)\endpmatrix\otimes \bfuno_s\in\RR^{sm},\\[2mm]
\label{notations2}
\bfomega(c) &=& \pmatrix{ccc} \omega_1(c) I_{m_1}\\ & \ddots \\ && \omega_\nu(c) I_{m_\nu}\endpmatrix\otimes I_s\in\RR^{sm\times sm}.
\end{eqnarray}
However, in order to derive a practical numerical method, the $s\nu$ Fourier coefficients (\ref{stepn}) (i.e.,  the block vector (\ref{stepn_v})) 
need to be approximated by suitable quadrature rules. In the case where in (\ref{fde1}) $\nu=1$, this can be efficiently done by using one  Gauss-Jacobi quadrature, defined at the zeros of the corresponding $k$th Jacobi polynomial: this results into a FHBVM$(k,s)$ method. In the following sections, we study the extension to the case $\nu>1$.

\section{Extending FHBVMs to the multi-order case}\label{nugt1}
The straightforward extension of FHBVMs to cope with problem (\ref{fde1}), consists in discretizing the integrals involved in the Fourier coefficients (\ref{stepn}) at the $n$th integration step, by means of a corresponding Gauss-Jacobi quadrature. In more details, for a suitable $k\ge s$, one approximates (\ref{stepn}) with the order $2k$ interpolatory quadrature rule defined at the zeros $c_\rho^i$,  $\rho=1,\dots,k$, of $P_k^i(c)$, with corresponding weights $\{b_\rho^i\}$: 
\begin{eqnarray}\nonumber
\gamma_{ij}(\sigma^n) &\approx& \gamma_{ij}^n~=~ \sum_{\rho=1}^k b_\rho^iP_j^i(c_\rho^i) f_i\left(\phi^n(c_\rho^ih_n)+h_n^\bfalfa I^\bfalfa P(c_\rho^i)\bfgamma^n\right),\\[1mm]
&& \qquad\qquad j=0,\dots,s-1,\qquad i=1,\dots,\nu, \label{stepnd}
\end{eqnarray}
having set (compare with (\ref{notations0}))
\begin{equation}\label{notations3}
\bfgamma^n = \pmatrix{c} \gamma_1^n\\ \vdots \\ \gamma_\nu^n\endpmatrix, \qquad \gamma_i^n=\pmatrix{c} \gamma_{i0}^n\\ \vdots \\ \gamma_{i,s-1}^n\endpmatrix,\qquad i=1,\dots,\nu.
\end{equation}
{\blue According to the analysis in Section~\ref{conval}, this} results into a FHBVM$(k,s)$ method {\blue (see also \cite{BBBI2024})}, having order at least $s$. In particular, the choice $k=s$, provides a collocation method. Moreover, a spectrally accurate solution can be achieved, by choosing $s\gg1$, in view of (\ref{spectacc}). This is exactly the strategy used  in  the two Matlab\cpr codes {\tt fhbvm} and {\tt fhbvm2}, implementing a FHBVM(22,20) and a FHBVM(22,22) method, respectively, designed for the case $\nu=1$. In the case $\nu>1$, however, one needs to evaluate the memory term, $\phi^n(ch_n)$, at the $k\nu$ abscissae
\begin{equation}\label{knuabs}
c_\rho^i, \qquad \rho=1,\dots,k, \qquad i=1,\dots,\nu.
\end{equation}
This is a quite expensive part of the computational cost since, as it can be guessed from (\ref{fijin}), all the information {\em ab initio} has to be computed for each abscissa in (\ref{knuabs}). Consequently, having $\nu$ different sets of abscissae is not a welcome news.

A different approach consists in looking for a set of suitable abscissae, $c_1,\dots,c_k$, such that the interpolatory quadrature rules defined over them are exact for polynomials of degree $2s-1$ (as it happens for the Gauss-Jacobi quadratures with $k=s$), for all the given weighting functions appearing in (\ref{stepn}). I.e.,\footnote{Hereafter, $\Pi_m$ will denote the vector space of polynomials of degree at most $m$.}  
\begin{equation}\label{mop1}
\forall g\in\Pi_{2s-1}~:~\int_0^1 \omega_i(c) g(c)\dd c ~=~ \sum_{\rho=1}^k b_\rho^i \,g(c_\rho), \quad i=1,\dots,\nu, 
\end{equation}
where, for the sake of brevity, we continue denoting by $\{b_\rho^i\}$ the weights of the new quadratures. 
This is equivalent to find a family of polynomials $\{\oP_k\}$, such that, for a suitable $q\in \NN$,
\begin{equation}\label{mop2}
\oP_k(c_\rho)=0,\quad \rho=1,\dots,k, \qquad \int_0^1 \omega_i(c)\oP_k(c) c^j\dd c = 0, \quad j=0,\dots,q-1, \quad i=1,\dots,\nu.
\end{equation}
This is a particular instance of a family of so-called {\em multiple orthogonal polynomials (MOPs)} \cite{VAC2001,CVA2005}, which we shall briefly sketch in the next section.
In particular, in the case of multiple Jacobi weights, as in our case, such abscissae are known as {\em \JP abscissae} \cite{LMVAVD2024}.

We observe that, in case (\ref{mop1})--(\ref{mop2}) hold true, then we can compute\,\footnote{Of course, the weights $\{b_\rho^i\}$ are not the same as those in (\ref{stepnd}) but, instead, aare those in (\ref{mop1}) .}
\begin{eqnarray}\nonumber
\gamma_{ij}(\sigma^n) &\approx& \gamma_{ij}^n~=~ \sum_{\rho=1}^k b_\rho^iP_j^i(c_\rho)f_i\left(\phi^n(c_\rho h_n)+h_n^\bfalfa I^\bfalfa P(c_\rho)\bfgamma^n\right),\\[1mm]
&& \qquad\qquad j=0,\dots,s-1,\qquad i=1,\dots,\nu, \label{stepnd1}
\end{eqnarray}
in place of (\ref{stepnd}). In such a way, we need to compute the memory term $\phi^n(ch_n)$ only at {\em one set} of $k$ abscissae, instead of $\nu$ sets.

\begin{rem}\label{iteg}
We observe that, when using (\ref{stepnd1}) in place of (\ref{stepnd}), also the polynomials $P_j^i(c)$, and the corresponding integrals $I^{\aa_i}P_j^i(c)$, $j=0,\dots,s-1$, $i=1,\dots,\nu$, are evaluated at only {\em one set} of abscissae, $c_1,\dots,c_k.$ Conversely, in the case of (\ref{stepnd}), they need to be evaluated at $\nu$ different sets of abscissae, and this impacts unfavorably  also in the implementation of the nonlinear iteration, as we shall sketch later in Section~\ref{newton}.
\end{rem}

\section{Multiple Orthogonal Polynomials (MOPs)}\label{mops}
To begin with, the aim is that of computing a set of abscissae $c_1,\dots,c_k$, satisying (\ref{mop2}), for a suitable $k$ to be determined. To this end, the following result holds true.

\begin{theo}\label{qth}
Assume (\ref{mop2}) holds true. Then the interpolatory quadratures with weights
$$b_\rho^i = \int_0^1 \omega_i(c)L_{\rho k}(c)\dd c, \qquad L_{\rho k}(c) = \prod_{\iota\ne \rho}\frac{c-c_\iota}{c_\rho-c_\iota}, \qquad \rho=1,\dots,k,\qquad i=1,\dots,\nu,$$
satisfy
$$\int_0^1\omega_i(c)g(c)\dd c  = \sum_{\rho=1}^k b_\rho^i g(c_\rho), \qquad i=1,\dots,\nu,\qquad \forall g\in\Pi_{k+q-1}.$$
\end{theo}
\proof See \cite[Theorem~1]{Bo1994}. \QED

\bigskip
Consequently, in order for this approach to reduce to the standard one \cite{BBBI2024,BGIV2025}, when $\nu=1$, hereafter we require:
\begin{equation}\label{req}
k = \nu q, \qquad k+q-1\ge 2s-1,
\end{equation}
so that all the quadratures $(c_\rho,b_\rho^i)$, $i=1,\dots,\nu$, have at least order $2s$. 
Therefore, the following result is proved.

\begin{cor}\label{qk}
The conditions (\ref{req}) are satisfied, provided that\,\footnote{Hereafter, $\lceil\cdot\rceil$ and $\lfloor\cdot\rfloor$ will denote the ceiling and floor functions, respectively.}
\begin{equation}\label{q}
q =  \left\lceil\frac{2s}{\nu+1}\right\rceil \qquad \Rightarrow\qquad k = \nu \left\lceil\frac{2s}{\nu+1}\right\rceil.
\end{equation}
\end{cor}
Table~\ref{knu} shows the values of $k$ for $s=22$, which is the value considered in the code {\tt fhbvm2} \cite{BGIV2025}.

\begin{table}[t]
\caption{Values of $k$ according to (\ref{q}) for increasing values of $\nu$, for $s=22$.}
\label{knu}
\centering

\smallskip
\begin{tabular}{|r|rrrrr|}
\hline\hline
$\nu$ & 1 & 2 & 3 & 4 & 5\\
\hline
$k$ & 22 & 30 & 33 & 36 & 40\\
\hline\hline
\end{tabular}
\end{table}

Consequently, the problem is that of finding the zeros of the polynomial $\oP_k(c)$ in (\ref{mop2}). A straightforward approach, to determine such a polynomial, consists in expressing it w.r.t. the power basis,
$$\oP_k(c) = \sum_{\iota=0}^k a_\iota c^{k-\iota}, \qquad a_0=1,$$
and impose the $k$ orthogonality conditions (\ref{mop2}) (with $k=\nu q$, according to (\ref{req})). These orthogonality conditions, by setting hereafter
\begin{equation}\label{scal}
(f(c),g(c))_i = \int_0^1\omega_i(c) f(c)g(c) \dd c, \qquad i=1,\dots,\nu,
\end{equation} 
read
$$\pmatrix{ccc}
(c^0,c^0)_1 & \dots & (c^0,c^{k-1})_1\\
\vdots          &          &\vdots\\
(c^{q-1},c^0)_1 & \dots & (c^{q-1},c^{k-1})_1\\
\vdots          &          &\vdots\\
(c^0,c^0)_\nu & \dots & (c^0,c^{k-1})_\nu\\
\vdots          &          &\vdots\\
(c^{q-1},c^0)_\nu & \dots & (c^{q-1},c^{k-1})_\nu
\endpmatrix\pmatrix{c} a_k\\ \vdots\\ a_1\endpmatrix =
-\pmatrix{c}
(c^0,c^k)_1\\ \vdots \\ (c^{q-1},c^k)_1\\
\vdots\\
(c^0,c^k)_\nu\\ \vdots \\ (c^{q-1},c^k)_\nu
\endpmatrix.
$$
Nevertheless, such linear system of $k$ equations has the pitfall of having a badly conditioned coefficient matrix and is, therefore, not practically useful to recover $\oP_k(c)$. A different approach consists in deriving this polynomial as the $k$th one in the sequence $\{\oP_j\}_{j\ge0}$ of monic polynomials satifying the recurrence \cite{MiSt2003}\,\footnote{The notation has been slightly modified, w.r.t. \cite{MiSt2003}, in order to have a more ``familiar'' enumeration of the entries of matrix $H_k$ in (\ref{Hk2}).}
\begin{equation}\label{bP}
\oP_j(c) = c\oP_{j-1}(c)- \sum_{i=\max(1,j-\nu)}^j  a_{ji} \oP_{i-1}(c), \quad j=1,2,\dots,k, \qquad
\oP_0(c)\equiv1.
\end{equation}
This can be cast in vector form as
\begin{equation}\label{Hk1}
H_k\pmatrix{c} \oP_0(c) \\ \oP_1(c) \\ \vdots \\ \vdots\\ \oP_{k-1}(c)\endpmatrix = 
c\pmatrix{c} \oP_0(c) \\ \oP_1(c) \\  \vdots \\ \vdots\\ \oP_{k-1}(c)\endpmatrix -
\pmatrix{c} 0 \\ \vdots \\ \vdots\\ 0\\ \oP_k(c)\endpmatrix,
\end{equation}
with $H_k\in\RR^{k\times k}$ a lower Hessenberg matrix,
\begin{equation}\label{Hk2}
H_k = \pmatrix{cccccc}
a_{11}         & 1\\
\vdots        & \ddots           & \ddots\\
a_{\nu+1,1} & \dots            & a_{\nu+1,\nu+1}        &1\\
                   &\ddots            &                          & \ddots                        &\ddots       \\
                   &                      & a_{k-1,k-1-\nu} &\dots                           &a_{k-1,k-1} &1\\
                   &                      &                          &a_{k,k-\nu}                 &\dots           &a_{kk}      
\endpmatrix.
\end{equation}
From (\ref{Hk1}) one deduces that $\oP_k(c_\rho)=0$, $\rho=1,\dots,k$, iff $c_\rho\in\sigma(H_k)$, with (right) eigenvectors respectively given by:
\begin{equation}\label{vro}
\pmatrix{c}
\oP_0(c_\rho)\\ \vdots \\ \oP_{k-1}(c_\rho)\endpmatrix,\qquad \rho=1,\dots,k. 
\end{equation}
The coefficients $\{a_{ji}\}$, in turn, if ~$j\in\{1,\dots,k\}$~ is such that
\begin{equation}\label{mu}
j=r\nu+\mu,\qquad \mbox{with}\qquad \mu\in\{0,\dots,\nu-1\},
\end{equation}
are computed by imposing that, by using the notation (\ref{scal}),\footnote{Of course, when $r-1<0$, it is meant that the corresponding orthogonality condition is void.} 
\begin{eqnarray*}
(\oP_j(c),c^\rho)_i&=&0,\qquad \rho=0,\dots,r, \qquad\qquad i=1,\dots,\mu,\\[2mm]
(\oP_j(c),c^\rho)_i&=&0,\qquad \rho=0,\dots,r-1, \,\qquad i=\mu+1,\dots,\nu,
\end{eqnarray*}
which is clearly equivalent to require:
\begin{eqnarray}\nonumber
(\oP_j(c),\oP_\rho(c))_i&=&0,\qquad \rho=0,\dots,r, \qquad\qquad i=1,\dots,\mu,\\[2mm]
(\oP_j(c),\oP_\rho(c))_i&=&0,\qquad \rho=0,\dots,r-1, \,\qquad i=\mu+1,\dots,\nu.\label{ortopi}
\end{eqnarray}
By setting, with reference to (\ref{mu}),
\begin{equation}\label{xij}
\xi(j) = \left\{\begin{array}{ccl}
\mu, &\mbox{if} &\mu>0,\\[2mm]
\nu,   &\mbox{if} & \mu=0,
\end{array}\right.
\end{equation}
this results into \cite{MiSt2003}: 
\begin{eqnarray}\nonumber
a_{ji} &=& \frac{\left(c\oP_{j-1}(c)-\sum_{\iota=\max(1,j-\nu)}^{i-1}a_{j\iota}\oP_{\iota-1}(c),\,\oP_{\left\lfloor \frac{i-1}\nu\right\rfloor}(c)\right)_{\xi(i)}}
{\left(\oP_{i-1}(c),\oP_{\left\lfloor \frac{i-1}\nu\right\rfloor}(c)\right)_{\xi(i)}},\\[2mm]
&&\qquad\qquad\qquad i=\max(1,j-\nu),\dots,j, \qquad j=1,\dots,k.\label{aji}
\end{eqnarray}
In fact, the following result holds true.

\begin{theo}\label{mopth}
By considering $k=q\nu$, according to (\ref{req}), the family of polynomials $\{\pi_1,\dots,\pi_k\}$,  defined by (\ref{bP})--(\ref{mu}) and (\ref{xij})--(\ref{aji}), satisfies the orthogonality conditions (\ref{ortopi}) (see (\ref{scal})). 
\end{theo}
\proof
Following an induction argument on the index $j$ varying on $\nu$ consecutive values, for $j=1$, one has 
$$\pi_1(c) = c-a_{11}, \qquad a_{11} = \frac{(c,1)_1}{(1,1)_1},$$
so that (\ref{ortopi}) holds true with $r=0$ and $\mu=1$. 
A straightforward induction argument, shows that for $j=2,\dots,\nu$, from (\ref{bP}) one obtains:
$$
(\pi_j(c),\pi_0(c))_i=0 \qquad \Rightarrow\qquad a_{ji} = \frac{\left(c\pi_{j-1}(c) - \sum_{\iota=1}^{i-1}a_{j\iota}\pi_{\iota-1}(c),\pi_0(c)\right)_i\qquad}{\left(\pi_{i-1}(c),\pi_0(c)\right)_i}, \qquad i=1,\dots,j,
$$
which matches (\ref{aji}), as required. Next, let us suppose (\ref{ortopi}) holds true up to $j=r\nu$, $r\ge1$, and prove for 
\begin{equation}\label{allmu}
j=r\nu+\mu, \qquad \mu=1,\dots,\nu. 
\end{equation}
Let us start with $\mu=1$. Clearly,
\begin{equation}\label{ximu1}
j=r\nu+1 \qquad \Rightarrow\qquad \xi(i) = 1,\dots,\nu,1, \qquad\mbox{as}\qquad i=j-\nu,\dots,j-1,j.
\end{equation}
Moreover, we know that (see (\ref{bP}))
$$i-1 \in\{ j-\nu-1,\dots,j-1\} = \{(r-1)\nu,\dots,r\nu\},$$
so that, by induction hypothesis, if ~$i-1=(r-1)\nu+\eta$, $\eta=0,\dots,\nu$,~ then the following orthogonality conditions hold true:\footnote{Clearly, if $r-2<0$, the corresponding ortogonality condition is void.}
\begin{eqnarray}\nonumber
(\oP_{j-\nu}(c),\oP_\rho(c))_\iota &=& 0, \qquad \rho=0,\dots,r-2, \qquad \iota=1,\dots,\nu,\qquad\qquad (\mbox{i.e.,} ~\eta=0),\\[5mm]
\nonumber
(\oP_{i-1}(c),\oP_\rho(c))_\iota &=& 0, \qquad \rho=0,\dots,r-1, \qquad \iota=1,\dots,\eta,\\[2mm]\nonumber
(\oP_{i-1}(c),\oP_\rho(c))_\iota &=& 0, \qquad \rho=0,\dots,r-2, \qquad \iota=\eta+1,\dots,\nu,
 \qquad \,\eta=1,\dots,\nu-1,\\[5mm]\label{rhoval}
(\oP_{j-1}(c),\oP_\rho(c))_\iota &=& 0, \qquad \rho=0,\dots,r-1, \qquad \iota=1,\dots,\nu,\qquad\qquad (\mbox{i.e.,} ~\eta=\nu).
 \end{eqnarray}
Next, according to (\ref{ortopi}),  let us impose that $\oP_j$ satisfies the following orthogonality conditions,
\begin{eqnarray*}
(\oP_j(c),\oP_\rho(c))_i&=&0, \qquad \rho=0,\dots,r-1,\qquad i=1,\dots,\nu,\\[2mm]
(\oP_j(c),\oP_r(c))_1&=&0,
\end{eqnarray*}
which, by virtue of (\ref{rhoval}), reduce to
\begin{eqnarray}\nonumber
(\oP_j(c),\oP_{r-1}(c))_i&=&0, \qquad i=1,\dots,\nu,\\[2mm]
(\oP_j(c),\oP_r(c))_1&=&0.\label{mu1}
\end{eqnarray}
Imposing such conditions in the given order, then gives, according to (\ref{ximu1}),
$$
a_{j,j-\nu+\eta} = \frac{\left(c\oP_{j-1}(c)-\sum_{\iota=j-\nu}^{j-\nu+\eta-1}a_{j\iota}\oP_{\iota-1}(c),\,\oP_{r-1}(c)\right)_{\eta+1}}
{\left(\oP_{j-\nu+\eta-1}(c),\oP_{r-1}(c)\right)_{\eta+1}}, \qquad \eta=0,\dots,\nu-1,
$$
and
$$
a_{jj} =  \frac{\left(c\oP_{j-1}(c)-\sum_{\iota=j-\nu}^{j-1}a_{j\iota}\oP_{\iota-1}(c),\,\oP_r(c)\right)_1}
{\left(\oP_{j-1}(c),\oP_r(c)\right)_1},
$$
respectively, which match (\ref{aji}) when $j=r\nu+1$. The subsequent values of $\mu$ in (\ref{allmu}) can be handled by means of similar arguments. 
In fact, for $\mu=2,\dots,\nu$, one needs to impose (compare with (\ref{mu1})):
\begin{eqnarray*}
(\oP_j(c),\oP_{r-1}(c))_i&=&0, \qquad i=\mu,\dots,\nu,\\[2mm]
(\oP_j(c),\oP_r(c))_i&=&0, \qquad i=1\dots,\mu.
\end{eqnarray*}
From these orthogonality conditions, by considering that, for ~$i=j-\nu,\dots,j,$ ~ $\xi(i)=\mu,\dots,\nu,1,\dots,\mu$,
 one eventually obtains:
$$
a_{j,j-\nu+\eta} = \frac{\left(c\oP_{j-1}(c)-\sum_{\iota=j-\nu}^{j-\nu+\eta-1}a_{j\iota}\oP_{\iota-1}(c),\,\oP_{r-1}(c)\right)_{\eta+\mu}}
{\left(\oP_{j-\nu+\eta-1}(c),\oP_{r-1}(c)\right)_{\eta+\mu}}, \qquad \eta=0,\dots,\nu-\mu,
$$
and
$$
a_{j,j-\nu+\eta} = \frac{\left(c\oP_{j-1}(c)-\sum_{\iota=j-\nu}^{j-\nu+\eta-1}a_{j\iota}\oP_{\iota-1}(c),\,\oP_r(c)\right)_{\eta-\nu+\mu}}
{\left(\oP_{j-\nu+\eta-1}(c),\oP_r(c)\right)_{\eta-\nu+\mu}}, \qquad \eta=\nu-\mu+1,\dots,\nu,
$$
which, again, are consistent with (\ref{aji}). Consequently, the statement follows.
\QED\bigskip

\begin{rem}
The proof of Theorem~\ref{mopth} has been reported for a twofold reason:
\begin{enumerate}
\item a formal proof, in the general case, was still lacking in the literature (usually, only the case $\nu=2$ is worked out);

\item the proof itself provides operative guidelines for efficiently deriving the entries of matrix $H_k$ in (\ref{Hk2}).
\end{enumerate}
\end{rem}

\subsection{Some computational remarks}
In order to derive a practical numerical method, one needs:
\begin{enumerate}
\item  at first, compute the integrands involved in (\ref{aji}) (see (\ref{scal})). By taking into account (\ref{req})--(\ref{q}), one has that the integrands are polynomials of degree at most
$$k-1+\left\lfloor\frac{k-1}{\nu}\right\rfloor = \nu q-1+ \left\lfloor\frac{\nu q-1}{\nu}\right\rfloor = \nu q +q-2 =(\nu+1)q-2= (\nu+1)\left\lceil\frac{2s}{\nu+1}\right\rceil-2.$$
Consequently, they can be computed exactly (up to round-off errors) by using a corresponding Gauss-Jacobi quadrature over 
\begin{equation}\label{vfi}
\varphi =\left\lceil \frac{\nu+1}2 \left\lceil\frac{2s}{\nu+1}\right\rceil-\frac{1}2 \right\rceil
\end{equation}
abscissae. For sake of completeness, in Table~\ref{fitab} we list a few values of $\varphi$, for the case $s=22$ considered in the code {\tt fhbvm2}, for increasing values of $\nu$;

\item secondly, matrix $H_k$ in (\ref{Hk1})--(\ref{Hk2}) must be balanced, in order to improve accuracy in numerically finding its eigenvalues. For this purpose, a recent strategy, proposed in \cite{LMVAVD2024}, consists in a diagonal scaling that symmetrizes the tridiagonal part of $H_k$. 
Additional computational aspects have been also considered in \cite{LMVD2025,VA2024}.

\end{enumerate}

\begin{table}[t]
\caption{Values of $\varphi$ in (\ref{vfi}) for increasing values of $\nu$, for $s=22$.}
\label{fitab}
\centering

\smallskip
\begin{tabular}{|r|rrrrr|}
\hline\hline
$\nu$ & 1 & 2 & 3 & 4 & 5\\
\hline
$\varphi$ & 22&  22&  22&  22&  24 \\
\hline\hline
\end{tabular}
\end{table}

\section{Solving the discrete problems}\label{newton}

As sketched above (see (\ref{stepnd1})), the discrete problem to be solved at the $n$th integration step consists in:
\begin{equation}\label{dispro}
\gamma_{ij}^n~=~ \sum_{\rho=1}^k b_\rho^iP_j^i(c_\rho)f_i\left(\phi^n(c_\rho h_n)+h_n^\bfalfa I^\bfalfa P(c_\rho)\bfgamma^n\right), \quad j=0,\dots,s-1,\quad i=1,\dots,\nu, 
\end{equation}
where the notations (\ref{notations0})--(\ref{notations1}) and (\ref{notations3}) have been used. By using notations (\ref{notations0})--(\ref{notations1}) and (\ref{notations3}) again, we can recast (\ref{dispro}) as:
\begin{equation}\label{dispro1}
\gamma_i^n = (\P_s^i)^\top\Omega_i\cdot f_i\left(\phi^n(\bfc h_n)+\h_n^\bfalfa \I_s^\bfalfa \bfgamma^n\right), \qquad i=1,\dots,\nu,
\end{equation}
where we have set
\begin{eqnarray}\nonumber
\bfc &=& \pmatrix{c} c_1\\ \vdots \\ c_k\endpmatrix\in\RR^k, \qquad 
\phi^n(\bfc h_n) = \pmatrix{c} \phi^n(c_1h_n)\\ \vdots \\ \phi^n(c_kh_n)\endpmatrix\in\RR^{km},\\[2mm] \label{notations4_0}
\h_n^\bfalfa &=& I_k\otimes h_n^\bfalfa\in\RR^{km\times km},
\end{eqnarray} 
additionally, 

\begin{eqnarray}\nonumber
\Omega_i &=& \pmatrix{ccc} b_1^i\\ &\ddots\\ && b_k^i\endpmatrix\otimes I_{m_i}\in\RR^{km_i\times km_i}, \\[2mm] 
\nonumber
\P_s^i &=&\pmatrix{ccc} P_0^i(c_1)& \dots &P_{s-1}^i(c_1)\\ \vdots & & \vdots\\ P_0^i(c_k)& \dots & P_{s-1}^i(c_k)\endpmatrix\otimes I_{m_i}\in\RR^{km_i\times sm_i},\\[2mm]
\I_s^\bfalfa &=& \pmatrix{c} I^\bfalfa P(c_1)\\ \vdots\\ I^\bfalfa P(c_k)\endpmatrix\in\RR^{km\times sm},\label{notations4}
\end{eqnarray}
and 
\begin{equation}\label{fisig}
f_i\left(\phi^n(\bfc h_n)+\h_n^\bfalfa \I_s^\bfalfa \bfgamma^n\right) = \pmatrix{c}
f_i\left(\phi^n(c_1 h_n)+h_n^\bfalfa I^\bfalfa P(c_1)\bfgamma^n\right)\\ \vdots\\
f_i\left(\phi^n(c_k h_n)+h_n^\bfalfa I^\bfalfa P(c_k)\bfgamma^n\right)\endpmatrix, \qquad i=1,\dots,\nu.
\end{equation}
By further introducing the notation
\begin{equation}\label{notations5}
\P_s = \pmatrix{ccc} \P_s^1\\ &\ddots \\ &&\P_s^\nu\endpmatrix\in\RR^{km\times sm},\qquad
\Omega = \pmatrix{ccc} \Omega_1\\ &\ddots\\ &&\Omega_\nu\endpmatrix\in\RR^{km\times km},
\end{equation}
and
\begin{equation}\label{fsig}
f\left(\phi^n(\bfc h_n)+\h_n^\bfalfa \I_s^\bfalfa \bfgamma^n\right) = \pmatrix{c}
f_1\left(\phi^n(\bfc h_n)+\h_n^\bfalfa \I_s^\bfalfa \bfgamma^n\right)\\
\vdots\\
f_\nu\left(\phi^n(\bfc h_n)+\h_n^\bfalfa \I_s^\bfalfa \bfgamma^n\right)\endpmatrix\in\RR^{km},
\end{equation}
one eventually obtains that (\ref{dispro1}) can be compactly written as:
\begin{equation}\label{dispro2}
\bfgamma^n = \P_s^\top\Omega \cdot f\left(\phi^n(\bfc h_n)+\h_n^\bfalfa \I_s^\bfalfa \bfgamma^n\right).
\end{equation}
From this formulation, the following result follows.

\begin{theo}\label{fix}
Assume $f$ is Lipschitz continuous. Then, for all sufficiently small $h_n>0$ problem (\ref{dispro2}) has a unique solution.
\end{theo}
\proof
In fact, from (\ref{notations0})--(\ref{notations1}) and (\ref{notations5}), it follows that, provided that, for any suitable matrix norm,
\begin{equation}\label{hL}
h_n^{\min_i\aa_i}L\|\P_s^\top\Omega\|\|\I_s^\bfalfa\|<1,
\end{equation} 
with $L$ the local Lipschitz constant, the iteration
\begin{equation}\label{fixit}
\bfgamma^{n,r+1} = \P_s^\top\Omega \cdot f\left(\phi^n(\bfc h_n)+\h_n^\bfalfa \I_s^\bfalfa \bfgamma^{n,r}\right), \qquad r=0,1,\dots,
\qquad \bfgamma^{n,0} =\bfzero,
\end{equation}
is a contraction.\,\QED\bigskip

\begin{rem} Even though the fixed-point iteration (\ref{fixit}) is straightforward and easy to implement, nevertheless, the condition (\ref{hL}) may be too restrictive,  if $L\gg1$: in such a case, a Newton-type iteration has to be used. This aspect will be studied in the next two subsections~\ref{nu_1} and \ref{nu_2}. \end{rem}

\subsection{The case $\nu=1$}\label{nu_1}
In the case $\nu=1$, i.e., when only one fractional derivative occurs, a very efficient Newton-type iteration can be considered, i.e., the so-called {\em blended iteration}, at first derived in the ODE case \cite{BM2002,BIT2011,book}, and then adapted for the FDE case \cite{BGI2024}. We briefly sketch it, for sake of completeness. In such a case, from (\ref{notations0})--(\ref{notations1}), (\ref{notations4_0})--(\ref{notations4}), and (\ref{notations5}), with $\nu=1$, considering that $m_1=m$, it follows that  (\ref{dispro2}) becomes
$$\gamma_1^n = (\P_s^1)^\top\Omega_1 \cdot f_1\left( \phi^n(\bfc h_n) + h_n^{\aa_1}\I_s^{\aa_1}\gamma_1^n\right).$$
Consequently, we need to solve the problem
$$G(\gamma) := \gamma-(\P_s^1)^\top\Omega_1 \cdot f_1\left( \phi^n(\bfc h_n) + h_n^{\aa_1}\I_s^{\aa_1}\gamma\right) = 0.$$
Some algebra shows that the simplified Newton iteration, for its solution, reads:\footnote{As is usual, $f_1'$ denotes the Jacobian matrix of $f_1$.}
\begin{eqnarray}\nonumber
\mbox{solve:}~\left[I_{sm}-h_n^{\aa_1}X_{11}\otimes f_1'(\phi^n(0))\right]\delta^j &=& -G(\gamma^j), \\[2mm]
\mbox{update:}~ \gamma^{j+1}&=&\gamma^j+\delta^j, \qquad j=0,1,\dots, \label{simpNewt}
\end{eqnarray}
with
\begin{eqnarray}\nonumber
X_{11} &=&  \pmatrix{ccc} 
P_0^1(c_1) & \dots &P_0^1(c_k)\\
\vdots & &\vdots\\
P_{s-1}^1(c_1) & \dots &P_{s-1}^1(c_k)
\endpmatrix
\pmatrix{ccc} b_1^1\\ &\ddots\\&&b_k^1\endpmatrix
\pmatrix{ccc} 
I^{\aa_1}P_0^1(c_1) & \dots &I^{\aa_1}P_{s-1}^1(c_1)\\
\vdots & &\vdots\\
I^{\aa_1}P_0^1(c_k) & \dots &I^{\aa_1}P_{s-1}^1(c_k)
\endpmatrix\\ \label{x11}
&&\in\RR^{s\times s},
\end{eqnarray}
where the initial condition ~$\gamma^0=0$~ can be conveniently used. The iteration (\ref{simpNewt})--(\ref{x11}), though straightforward, requires to factor the matrix 
\begin{equation}\label{smmat}
\left[I_{sm}-h_n^{\aa_1}X_{11}\otimes f_1'(\phi^n(0))\right]~\in~\RR^{sm\times sm}.
\end{equation}
The associated {\em blended iteration}, formally obtained by combining (\ref{simpNewt}) with the equivalent formulation
$$\rho_{11}\left[X_{11}^{-1}\otimes I_m-h_n^{\aa_1} I_s\otimes f_1'(\phi^n(0))\right]\delta_j = -\rho_{11}(X_{11}^{-1}\otimes I_m) G(\gamma^j),$$
with weights ~$\theta^n$~ and ~$I_{sm}- \theta^n$,~ respectively, where
$$\theta^n = I_s\otimes \left[I_m-\rho_{11}h_n^{\aa_1}f_1'(\phi^n(0))\right]^{-1}$$
and $\rho_{11}>0$ is a suitable parameter depending on the spectrum of matrix $X_{11}$ \cite{BGI2024}, reads, as follows:
\begin{eqnarray}\nonumber
\mbox{set:\qquad~} \eta_j &=& -G(\gamma^j)\\[2mm] \nonumber
\hat\eta_j &=& \rho_{11}(X_{11}^{-1}\otimes I_m)\eta_j\\[2mm] \label{blend}
 \gamma^{j+1} &=&\gamma^j-\theta^n\left[ \hat\eta_j + \theta^n\left(\eta_j-\hat\eta_j\right)\right], \qquad j=0,1,\dots.
\end{eqnarray}
Also in this case, one conveniently starts at $\gamma^0=0$. Clearly, this iteration only requires to factor the matrix 
$$I_m-\rho_{11}h_n^{\aa_1}f_1'(\phi^n(0)) ~\in~\RR^{m\times m},$$
instead of (\ref{smmat}). Unfortunately, this very effective approach, implemented in the Matlab\cpr codes {\tt fhbvm} \cite{BGI2024} and {\tt fhbbm2} \cite{BGIV2025}, does not extend to the case $\nu>1$.

\subsection{The case $\nu>1$}\label{nu_2}
To cope with this case, we need to rely on the straight simplified Newton iteration for solving (\ref{dispro2}), i.e.,
$$G(\bfgamma) := \bfgamma - \P_s^\top\Omega\cdot f\left(\phi^n(\bfc h_n) +\bfh_n^\bfalfa \I_s^\bfalfa\bfgamma\right) = \bfzero,$$
which now reads:
\begin{eqnarray}\nonumber
\mbox{solve:}~\left[I_{sm}- \P_s^\top\Omega \cdot (\,I_k\otimes f'(\phi^n(0))\,) \cdot \bfh_n^\bfalfa\I_s\right]\bfdelta^j &=& -G(\bfgamma^j), \\[2mm]
\mbox{update:}~ \bfgamma^{j+1}&=&\bfgamma^j+\bfdelta^j, \qquad j=0,1,\dots,\label{simpNewt1}
\end{eqnarray}
where, again, ~$\bfgamma^0=\bfzero$~ is conveniently used. For sake of clarity, let us explain in detail the (block) structure of the matrix
$$\P_s^\top\Omega \cdot (\,I_k\otimes f'(\phi^n(0))\,)\cdot \bfh_n^\bfalfa\I_s ~\equiv~ \pmatrix{ccc}
h_n^{\aa_1}X_{11}\otimes f_{11}^n & \dots &h_n^{\aa_\nu}X_{1\nu}\otimes f_{1\nu}^n\\
\vdots & &\vdots\\
h_n^{\aa_1}X_{\nu 1}\otimes f_{\nu 1}^n & \dots &h_n^{\aa_\nu}X_{\nu\nu}\otimes f_{\nu\nu}^n
\endpmatrix\in\RR^{sm\times sm},$$ 
where, with reference to (\ref{fde1}) and (\ref{fsig}), for all ~$i,j=1,\dots,\nu$:
$$f_{ij}^n = \frac{\partial}{\partial y_j}f_i(\phi^n(0))\in\RR^{m_i\times m_j},$$
and, generalizing (\ref{x11}),
\begin{eqnarray}\nonumber
X_{ij} &=& \pmatrix{ccc} 
P_0^i(c_1) & \dots &P_0^i(c_k)\\
\vdots & &\vdots\\
P_{s-1}^i(c_1) & \dots &P_{s-1}^i(c_k)
\endpmatrix
\pmatrix{ccc} b_1^i\\ &\ddots\\&&b_k^i\endpmatrix
\pmatrix{ccc} 
I^{\aa_j}P_0^j(c_1) & \dots &I^{\aa_j}P_{s-1}^j(c_1)\\
\vdots & &\vdots\\
I^{\aa_j}P_0^j(c_k) & \dots &I^{\aa_j}P_{s-1}^j(c_k)
\endpmatrix\\
&&\in\RR^{s\times s}.\label{Xijd}
\end{eqnarray}
Clearly, in this case the $sm\times sm$ coefficient matrix in (\ref{simpNewt1}), i.e.,
\begin{equation}\label{smxsm} 
 I_{sm}-\pmatrix{ccc}
h_n^{\aa_1}X_{11}\otimes f_{11}^n & \dots &h_n^{\aa_\nu}X_{1\nu}\otimes f_{1\nu}^n\\
\vdots & &\vdots\\
h_n^{\aa_1}X_{\nu 1}\otimes f_{\nu 1}^n & \dots &h_n^{\aa_\nu}X_{\nu\nu}\otimes f_{\nu\nu}^n
\endpmatrix\in\RR^{sm\times sm},
\end{equation} 
needs to be factored and, therefore, the computational cost is higher than in the case $\nu=1$, where the more efficient {\em blended iteration} (\ref{blend}) can be used.

\begin{rem}\label{GJquad}
In case one uses the approach (\ref{stepnd}), instead of (\ref{stepnd1}), for approximating the Fourier coefficients, the simplified Newton iteration would require to factor the matrix
$$I_{sm}- \pmatrix{ccc}
h_n^{\aa_1}X_{11}^1\otimes f_{11}^n & \dots &h_n^{\aa_\nu}X_{1\nu}^1\otimes f_{1\nu}^n\\
\vdots & &\vdots\\
h_n^{\aa_1}X_{\nu 1}^\nu\otimes f_{\nu 1}^n & \dots &h_n^{\aa_\nu}X_{\nu\nu}^\nu\otimes f_{\nu\nu}^n
\endpmatrix\in\RR^{sm\times sm},$$ 
in place of (\ref{smxsm}), with
\begin{eqnarray*}
X_{ij}^i &=& \pmatrix{ccc} 
P_0^i(c_1^i) & \dots &P_0^i(c_k^i)\\
\vdots & &\vdots\\
P_{s-1}^i(c_1^i) & \dots &P_{s-1}^i(c_k^i)
\endpmatrix
\pmatrix{ccc} b_1^i\\ &\ddots\\&&b_k^i\endpmatrix
\pmatrix{ccc} 
I^{\aa_j}P_0^j(c_1^i) & \dots &I^{\aa_j}P_{s-1}^j(c_1^i)\\
\vdots & &\vdots\\
I^{\aa_j}P_0^j(c_k^i) & \dots &I^{\aa_j}P_{s-1}^j(c_k^i)
\endpmatrix\\
&&\in\RR^{s\times s}.
\end{eqnarray*}
Consequently, for all \,$i=1,\dots,\nu$,\, the \,$i$th block row depends on the corrresponding set of Jacobi abscissae (i.e., the zeros of $P_k^i(c)$).
\end{rem}

{\blue

\section{Analysis of the methods}\label{analmet}

In this section we sketch a linear stability analysis of FHBVMs and  study their convergence order. We also report a few numerical tests to assess the theoretical findings. 

\subsection{Linear stability analysis}\label{linal}
We consider problem (\ref{fde1}) in the linear and autonomous case, assuming, without loosing generality, that $\nu=m$:
\begin{equation}\label{p1}
y_i^{(\aa_i)}(t) = \sum_{j=1}^m a_{ij} y_j(t), \qquad t>0, \qquad y_i(0)=y_{i0}\in\RR,\qquad 0<\aa_i<1, \qquad i=1,\dots,m.
\end{equation}
We also consider the corresponding shifted orthonormal Jacobi polynomial bases (\ref{ortoi}). Then, fixing a suitable $s\ge1$, according to (\ref{sigalfa})--(\ref{spectacc}), 
a FHBVM, used with timestep $t$, basically determines approximations:
\begin{equation}\label{ssi}
\sigma_i^{(\aa_i)}(ct) = \sum_{\nu=0}^{s-1} P_\nu^i(c)\gamma_{i\nu}, 
\qquad \sigma_i(ct) = y_{i0} + t^{\aa_i}\sum_{\nu=0}^{s-1} I^{\aa_i} P_\nu^i(c)\gamma_{i\nu},\qquad c\in[0,1], 
\end{equation}
with the coefficients $\{\gamma_{i\nu}\}$ determined such that  
\begin{equation}\label{o2}
\int_0^1 \omega_i(c)P_\iota^i(c)\left[ \sigma_i^{(\aa_i)}(ct) - \sum_{j=1}^m a_{ij} \sigma_j(ct)\right]\dd c = 0, \qquad \iota = 0,\dots,s-1, \quad i=1,\dots,m.
\end{equation}
From (\ref{ortoi}) and (\ref{ssi}), one then obtains that, recalling that $\delta_{\iota0}$ denotes the Kronecker symbol,
\begin{equation}\label{o3}
\gamma_{i\iota} - \sum_{j=1}^m a_{ij} \left[ \delta_{\iota0}y_{j0} +t^{\aa_j}\sum_{\nu=0}^{s-1}\left(\int_0^1\omega_i(c)P_\iota^i(c) I^{\aa_j}P_\nu^j(c)\dd c\right) \gamma_{j\nu}\right]=0,\quad \iota=0,\dots,s-1, \quad i=1,\dots,m.
\end{equation}
More compactly, by setting
\begin{equation}\label{fi}
x_{ij}^{\iota\nu} = \int_0^1\omega_i(c)P_\iota^i(c) I^{\aa_j}P_\nu^j(c)\dd c, \qquad i,j=1,\dots,m, \quad \iota,\nu=0,\dots,s-1,
\end{equation}
and
\begin{equation}\label{psi}
\psi_i = \sum_{j=1}^m a_{ij}y_{j0}, \qquad i=1,\dots,m,
\end{equation}
one derives that (\ref{o3}) is equivalent to:
\begin{equation}\label{o4}
\gamma_{i\iota} - \sum_{j=1}^m a_{ij} t^{\aa_j} \sum_{\nu=0}^{s-1} x_{ij}^{\iota\nu}\gamma_\nu^j = \delta_{\iota0}\psi_i,\qquad \iota = 0,\dots,s-1,
\quad  i=1,\dots,m.
\end{equation}
Setting $e_1\in\RR^s$ the first unit vector, the vectors
\begin{equation}\label{gi}
\gamma^i = \pmatrix{c} \gamma_{i0}\\ \vdots \\ \gamma_{i,s-1}\endpmatrix\in\RR^s, \quad i=1,\dots,m,
\end{equation}
the matrices 
\begin{equation}\label{Xij}
X_{ij} = \pmatrix{ccc} x_{ij}^{0,0} &\dots &x_{ij}^{0,s-1}\\
\vdots& &\vdots\\  x_{ij}^{s-1,0} &\dots &x_{ij}^{s-1,s-1}\endpmatrix, \quad B_{ij} = a_{ij}X_{ij}~\in~\RR^{s\times s}, \qquad i,j=1,\dots,m,
\end{equation}
and recalling $I_s\in\RR^{s\times s}$ the identity matrix, one has that (\ref{o4}) can be cast in vector form as
\begin{equation}\label{o5}
\pmatrix{cccc}
I_s-B_{11}t^{\aa_1} & -B_{12}t^{\aa_2}    &\dots &-B_{1m}t^{\aa_m}\\
-B_{21}t^{\aa_1}    &I_s-B_{22}t^{\aa_2} &\dots & -B_{2m}t^{\aa_m}\\
\vdots                    &                               &\ddots &\vdots\\
-B_{m1}t^{\aa_1}   & -B_{m2}t^{\aa_2} &\dots  & I_s-B_{mm}t^{\aa_m}\endpmatrix 
\pmatrix{c}\gamma^1\\ \gamma^2\\ \vdots\\ \gamma^m\endpmatrix
= \pmatrix{c} e_1\otimes\psi_1\\ e_1\otimes \psi_2\\ \vdots\\ e_1\otimes \psi_m\endpmatrix.
\end{equation}
\begin{rem}\label{xijrem}
Discretizing the integrals in (\ref{fi}) with a quadrature over $k$ points (e.g., $k$ Gauss-Pi\~{n}eiro abscissae corresponding to the given $\{\aa_i\}$) results into a FHBVM$(k,s)$ method. In such a case, the matrices $X_{ij}$ defined in (\ref{Xij}) do coincide with the matrices $X_{ij}$ defined in (\ref{Xijd}). Conversely, if a different Gauss-Jacobi quadrature is used for each equation, then $X_{ij}$ reduce to $X_{ij}^i$, as defined in Remark~\ref{GJquad}.
\end{rem}

Considering that the functions $\sigma_i(ct)$ in (\ref{ssi}) can be rewritten as
$$\sigma_i(ct) = y_{i0} + \sum_{\nu=0}^{s-1} I^{\aa_i} P_\nu^i(c)\left(t^{\aa_i}\gamma_{i\nu}\right),\qquad c\in[0,1],
 \qquad i=1,\dots,m,$$
it makes sense to reformulate problem (\ref{o5}) in the equivalent form:
\begin{equation}\label{o6}
M_s^\bfalfa (t)\pmatrix{c}t^{\aa_1}\gamma^1\\ t^{\aa_2} \gamma^2\\ \vdots\\ t^{\aa_m} \gamma^m\endpmatrix
= \pmatrix{c} e_1\otimes\psi_1\\ e_1\otimes \psi_2\\ \vdots\\ e_1\otimes \psi_m\endpmatrix,
\end{equation}
with ~$\bfalfa=(\aa_1,\dots,\aa_m)$~ and
\begin{equation}\label{M}
M_s^\bfalfa(t) = \pmatrix{cccc}
t^{-\aa_1}I_s-B_{11} & -B_{12}    &\dots &-B_{1m}\\
-B_{21}    &t^{-\aa_2}I_s-B_{22} &\dots & -B_{2m}\\
\vdots                    &                               &\ddots &\vdots\\
-B_{m1}   & -B_{m2} &\dots  & t^{-\aa_m}I_s-B_{mm}\endpmatrix.
\end{equation}
Consequently,
\begin{equation}\label{o7}
\pmatrix{c}t^{\aa_1}\gamma^1\\ t^{\aa_2} \gamma^2\\ \vdots\\ t^{\aa_m} \gamma^m\endpmatrix
= M_s^\bfalfa (t)^{-1}\pmatrix{c} e_1\otimes\psi_1\\ e_1\otimes \psi_2\\ \vdots\\ e_1\otimes \psi_m\endpmatrix,
\end{equation}
Considering that the r.h.s. in (\ref{o6}) does not depend on $t$, from (\ref{o7}) one deduces that a FHBVM$(k,s)$ method can be called {\em $A$-stable} if,  for all matrices ~$A=(a_{ij})\in\RR^{m\times m}$~ and ~$\bfalfa=(\aa_1,\dots,\aa_m)$~ such that the solution of problem (\ref{p1}) is asymptotically stable, one has that matrix
$$M_s^\bfalfa(t)^{-1}D_s, \qquad {\rm with}\qquad D_s = \diag( u\otimes e_1 ), \qquad u = (1,\dots,1)^\top\in\RR^m,$$ has a bounded norm for all $t>0$ (the reason for introducing matrix $D_s$ is due to the fact that, because of the structure of the r.h.s.,  only the columns $1,s+1,2s+1,\dots,(m-1)s+1$ of the inverse are actually involved). In this case, in fact, it follows that:
$$\exists\, K>0 ~ s.t. ~  \forall \, i=1,\dots,m, ~ t>0, \mbox{~and~} c\in[0,1]: ~ |\sigma_i(ct)|\le K.$$

As is clear, in general it is not possible to directly assess the $A$-stability of FHBVMs (as it has been done in the case of only one fractional order in \cite{BGIV2025}). Instead, we provide an example showing that FHBVM$(30,s)$, $s=1,\dots,22$, behaves correctly when solving the linear problem
\begin{equation}\label{stab}
y_1^{(0.5)}(t) = 0.1 y_1(t)+y_2(t),\quad y_2^{(0.7)}(t) = -y_1(t) + 0.1 y_2(t), \quad t>0, \qquad y_1(0)=y_2(0)=1,
\end{equation}
whose solution, according to the criteria given, e.g., in \cite{Deng} (see also \cite{BGK}), is asymptotically stable.\footnote{\blue In fact, in the multi-order case, the much simpler criterion in \cite{Matignon}, valid in the single-order case, cannot be used.} In more detail, in Figure~\ref{stabfig} we plot the function 
\begin{equation}\label{Sigmas}
\Sigma_s(t) := \| (M_s^{(0.5,0.7)}(t))^{-1}D_s\|, \qquad t\in[10^{-3},10^4], \qquad s=1,\dots,22.
\end{equation}
\begin{figure} 
\vspace{-2.5cm}
\centering
\includegraphics[width=16cm]{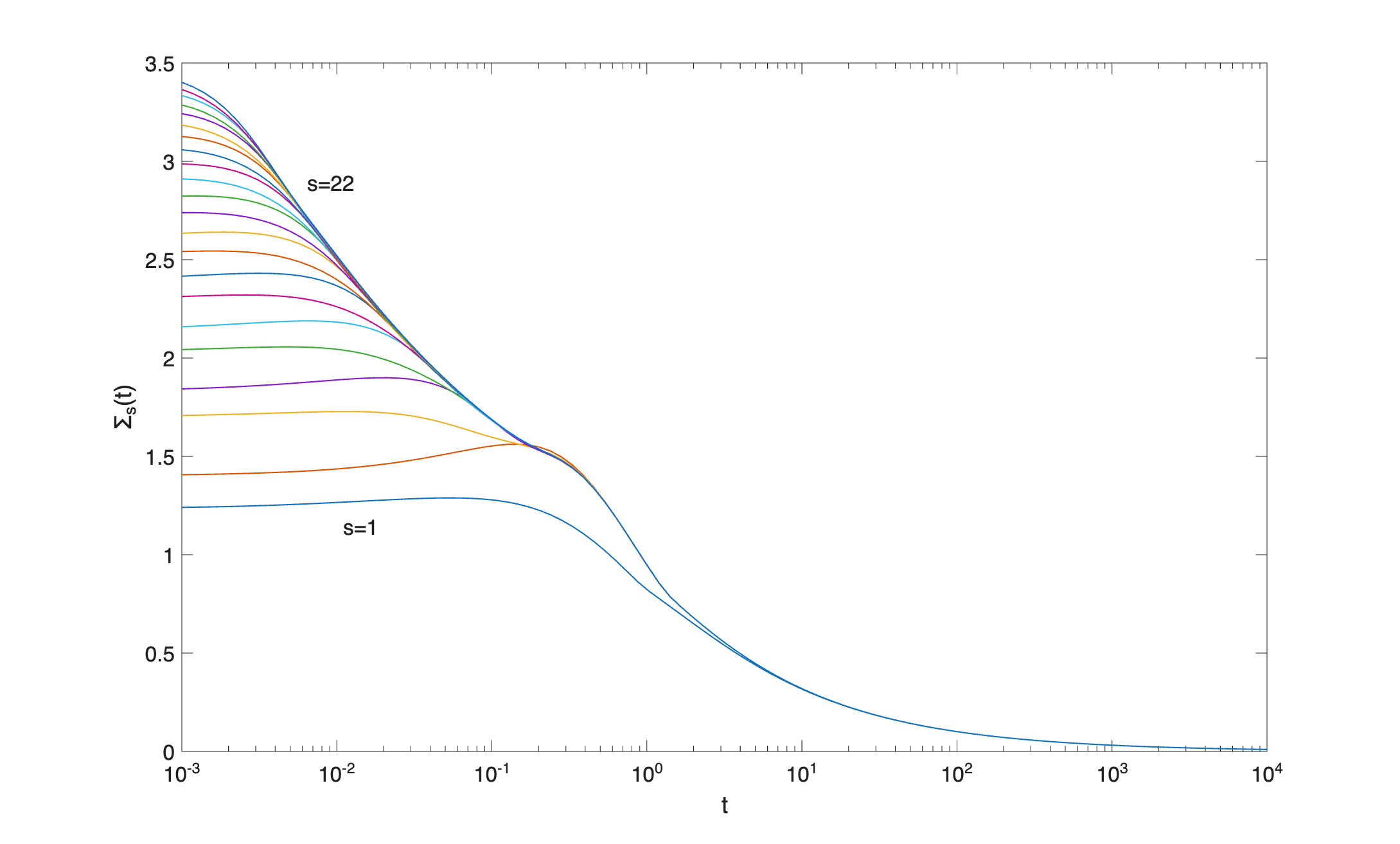} 
\caption{\blue Plot of the function $\Sigma_s(t)$ defined in (\ref{Sigmas}), $s=1,\dots,22$.}
\label{stabfig}
\end{figure}
From the plots, one deduces that all methods behave in the correct way.

\subsection{Convergence analysis}\label{conval}
In the sequel (see also \cite{BGIV2025}), we shall consider the following meshes:
\begin{itemize}
\item either a {\em graded mesh}, where
\begin{equation}\label{graded}
h_n = h_1 r^{n-1}, \qquad n=1,\dots,N, \qquad \mbox{with}\qquad r>1\mbox{\quad s.t.\quad} h_1\frac{r^N-1}{r-1}=T,
\end{equation}
or
\item a {\em uniform mesh}
\begin{equation}\label{uniform}
h_n = h \equiv \frac{T}N, \qquad n=1,\dots,N,
\end{equation}
or
\item a {\em mixed-mesh} where, for a suitably chosen \,$\rho\in\NN$,\, $\rho\ge1$,\, and setting
\begin{equation}\label{rmixed}
r \,:=\, \frac{\max\{2,\rho\}}{\max\{2,\rho\}-1},
\end{equation}
one has, for a suitable \,$\mu\in\NN$,\, $\mu\ge1$,\, and denoting \,$M:=N-\mu+\rho$\,:
\begin{eqnarray}\nonumber
h_n&=& h_1 r^{n-1}, ~~\qquad n=1,\dots,\mu,\qquad \mbox{with}\qquad h_1\frac{r^\mu-1}{r-1}=\rho h,\\
h_n&=&h~\equiv~\frac{T}{M},\qquad n=\mu+1,\dots,N.\label{mixed} 
\end{eqnarray}
\end{itemize}
Of course, (\ref{mixed}) reduces to either (\ref{uniform}), when $\mu=\rho=1$, or to (\ref{graded}), when $\mu=N>1$. Moreover, as is clear \cite{BGIV2025}:
\begin{itemize}
\item (\ref{graded}) is appropriate when the vector field is nonsmooth at $t=0$;
\item (\ref{uniform}) is appropriate when the vector field is suitably regular everywhere;
\item (\ref{mixed}) is appropriate when the vector field is nonsmooth at $t=0$, but the solution becomes eventually periodic.
\end{itemize}
In the sequel we shall report a unified convergence analysis for a HBVM$(k,s)$ method assuming, without loss of generality, that ~$T\le1$~ in (\ref{fde1}), in order to simplify the arguments.\footnote{\blue Differently, a time scaling can be done.}   By taking into account (\ref{fde1i}), (\ref{sigalfa}), (\ref{sigin}), (\ref{stepn}),  (\ref{stepnd}), and (\ref{stepnd1}), and considering that, for all $n=1,\dots,N$,
$$f_i(\sigma^n(ch_n)) = \sum_{j\ge0} P_j^i(c)\gamma_{ij}(\sigma^n), \qquad c\in[0,1], \qquad i=1,\dots,\nu,$$
 one obtains that the discrete solution exactly solves the sequence of perturbed problems
\begin{eqnarray}\nonumber
\sigma_{in}^{(\aa_i)}(ch_n) &=& f_i(\sigma^n(ch_n)) - \psi_{in}(ch_n), \qquad c\in[0,1],\qquad n=1,\dots,N,\\ \nonumber 
\sigma_{i1}^{(\iota)}(0) &=& y_{i0}^\iota\in\RR^{m_i},\qquad \iota=0,\dots,\ell-1, \qquad i=1,\dots,\nu,\\[2mm]
\sigma^n &=& (\sigma_{1n}^\top,\dots,\sigma_{\nu n}^\top)^\top\in\RR^m,\label{sigalfa1}
\end{eqnarray}
in place of (setting $y_{in}(ch_n) = y_i(t_{n-1}+ch_n)$, ~$c\in[0,1]$): 
\begin{eqnarray}\nonumber
y_{in}^{(\aa_i)}(ch_n) &=& f_i(y^n(ch_n)), \qquad c\in[0,1],\qquad n=1,\dots,N,\\ \nonumber 
y_{i1}^{(\iota)}(0) &=& y_{i0}^\iota\in\RR^{m_i},\qquad \iota=0,\dots,\ell-1, \qquad i=1,\dots,\nu,\\[2mm]
y^n &=& (y_{1n}^\top,\dots,y_{\nu n}^\top)^\top\in\RR^m.\label{yalfa1}
\end{eqnarray}
In (\ref{sigalfa1}), we have denoted
\begin{equation}\label{Psiin}
\psi_{in}(ch_n) ~=~ \sum_{j\ge s} P_j^i(c)\gamma_{ij}(\sigma^n)\,+\, \sum_{j=0}^{s-1} P_j^i(c)\left[\gamma_{ij}(\sigma^n)-\gamma_{ij}^n\right],
\end{equation}
the perturbation term, which is given by the sum of the truncation error in the expansion of $f(\sigma^n(ch_n))$ and of the quadrature error in the approximation of the Fourier coefficients  $\gamma_{ij}(\sigma^n)$ (see (\ref{sigin})).
In this respect, hereafter we shall assume that $k$ is large enough so that the quadrature implemented in the FHBVM$(k,s)$ method has order at least $2s$. In such a case, following arguments similar to those in \cite[Section\,3]{BBBI2024}, it can be shown that, under regularity assumptions on the functions $f_i$, one has that\,\footnote{\blue Hereafter, $\|\cdot\|$ will denote the infinity norm.} 
\begin{equation}\label{Psiin1}
\Psi_{in} := \max_{c\in[0,1]}\|\psi_{in}(ch_n)\|,
\end{equation}
satisfies
\begin{itemize}
\item for (\ref{graded}):
$$\Psi_{in}~\le~\left\{\begin{array}{ccl} O(h_1^{\aa_i}), &~&\mbox{if}~ n=1,\\[2mm] O(h_n^s),&&\mbox{otherwise};\end{array}\right.$$
\item for (\ref{uniform}): \hspace{3.5cm}$\Psi_{in}~\le~ O(h^s)$;
\item for (\ref{mixed}):
$$\Psi_{in}~\le~\left\{\begin{array}{ccl} O(h_1^{\aa_i}), &~&\mbox{if}~ n=1,\\[2mm] O(h_n^s),&&\mbox{if}~n=2,\dots,\mu,\\[2mm]
O(h^s),&&\mbox{otherwise.}\end{array}\right.$$
\end{itemize}
Setting ~$\Psi^n=\max_i\Psi_{in}$,
$$e_{in}(ch_n) = \|y_{in}(ch_n)-\sigma_{in}(ch_n)\|, \qquad e^n(ch_n) = \max_i e_{in}(ch_n)\qquad c\in[0,1], \qquad n=1,\dots,N,$$ 
the piecewise functions
$$ e(t_{n-1}+ch_n)\equiv e^n(ch_n),\qquad \Psi(t_{n-1}+ch_n) \equiv \Psi^n, \qquad c\in[0,1], \qquad n=1,\dots,N,$$
$\aa = \min_i\aa_i$, $L$ the maximum of the Lipschitz constants of the $f_i$, and considering that $\Gamma(x)^{-1}<1.2$, for $x>0$, one eventually arrives at the inequality:
$$e(t)~\le~ 1.2\int_0^t (t-\tau)^{\aa-1}\left[ L e(\tau) + \Psi(\tau)\right]\dd\tau, \qquad t\in[0,T].$$
We observe that, because of the choices (\ref{graded})--(\ref{mixed}) one deduces that:
$$
\forall\,2,\dots,N~:~ h_{n-1} \le h_n.
$$
Consequently, we can assume that $\Psi(t)$ is non decreasing and therefore, by using the fractional version of the Gronwall Lemma \cite[Corollary\,2]{YGD2006} (see also \cite{webb2021}), one deduces that:
\begin{equation}\label{et}
e(t)~\le~ K\int_0^t (t-\tau)^{\aa-1}\Psi(\tau)\dd\tau,
\end{equation}
with ~$K = 1.2\,E_\aa(1.2\Gamma(\aa)LT^\aa)$,\, being ~$E_\aa(z) = \sum_{j\ge0}\frac{z^j}{\Gamma(\aa j+1)}$~ the one-parameter Mittag-Leffler function.
Furthermore, considering that\,\footnote{\blue Recall that ~$h_n = t_n-t_{n-1}$~ and, moreover, ~$(t-t_{n-1})\le T\le 1$.}
$$\forall\, n=1,2,\dots,N, \mbox{\quad and\quad} t\ge t_n~:~ \int_{t_{n-1}}^{t_n}(t-\tau)^{\aa-1}\dd\tau ~\le~ \frac{ h_n^\taa}{\taa} ,$$ 
with
\begin{equation}\label{talfa}
\taa = \min\{\aa,1\},
\end{equation}
and taking into account that $\Psi(t)$ is piecewise constant, one obtains:
$$K\int_0^t (t-\tau)^{\aa-1}\Psi(\tau) \dd \tau ~\le~ \frac{K}\taa \sum_{k=1}^N h_k^\taa\Psi^k.$$
Moreover, by considering that for the graded mesh (\ref{graded}) and the graded part in (\ref{mixed})
$$\sum_{k=1}^n h_k^{s+\taa} = h_n^{s+\taa}\sum_{k=0}^{n-1} r^{-k(s+\taa)} < \frac{h_n^{s+\taa}}{1-r^{-s-\taa}},$$
whereas for the uniform mesh $nh^{s+\taa} \le h^{s+\taa-1}$, from (\ref{Psiin1}) one eventually obtains that:
\begin{equation}\label{finerr}
\max_{0\le t\le T}|e(t)|~\le~\left\{\begin{array}{ccl}
O(h_1^{\aa+\taa}) + O(h_N^{s+\taa}),&&\mbox{for}\quad (\ref{graded}),\\[2mm]
O(h^{s+\taa-1}),&&\mbox{for}\quad (\ref{uniform}),\\[2mm]
O(h_1^{\aa+\taa}) + O(h_\mu^{s+\taa})+O(h^{s+\taa-1}),&&\mbox{for}\quad (\ref{mixed}).\\[2mm]
\end{array}\right.
\end{equation}
In Section~\ref{numtest} we report a few numerical tests, to assess (\ref{talfa})--(\ref{finerr}) in the case (\ref{graded}).
\begin{rem}
From (\ref{finerr}) one concludes that  one can achieve arbitrarily high-order methods (though, suitably choosing the initial timestep, in the case of the graded or mixed meshes): as matter of fact, if $s\gg1$, a {\em spectrally accurate} solution can be gained, as it happens in the ODE case \cite{ABI2020}.
\end{rem}
}

\section{The code \fhbvmnew}\label{code}

{\blue
The Matlab\cpr code \fhbvmnew\, implements a  FHBVM$(k,22)$ method, using the approach described in the previous sections, in the case $\nu\in\{1,2\}$ in (\ref{fde1}). The code may  use one of the meshes (\ref{graded})--(\ref{mixed}) described in Section~\ref{conval}.} Concerning the evaluation of the abscissae and the weights for the quadrature(s):
\begin{itemize}
\item when $\nu=1$, {\blue $k=22$  and} the Jacobi abscissae and weights are computed by adapting the companion Matlab\cpr code of \cite{Gautschi2004};

\item when $\nu=2$, {\blue $k=30$ and} the \JP abscissae and weights are computed by adapting the companion Matlab\cpr code of \cite{LMVAVD2024}.
\end{itemize}
 {\blue 
 \begin{rem}
 We emphasize that the code \fhbvmnew\, currently handles only the cases $\nu=1,2$, because the software in \cite{LMVAVD2024} only implements the case of two weighting functions. On the other hand, the extension of such software is not trivial, and is beyond the scope of this paper. Consequently, when an extended version of the code in \cite{LMVAVD2024} will be available, we shall adapt the code \fhbvmnew\, accordingly, using the arguments explained in Sections~\ref{mops} and \ref{newton}.
 \end{rem}
 }

\no Concerning the nonlinear iteration:
\begin{itemize}
\item when $\nu=1$, the more efficient {\em blended iteration} (\ref{blend}) is used (as an alternative to the fixed-point iteration (\ref{fixit}), used whenever possible);

\item when $\nu=2$, the straight simplified-Newton iteration (\ref{simpNewt1}) is used (as an alternative to the fixed-point iteration (\ref{fixit}), used whenever possible).
\end{itemize}

\no All the remaining implementation details are similar to those of the codes {\tt fhbvm} and {\tt fhbvm2}, as explained in \cite{BGI2024} and \cite{BGIV2025}, respectively.  The new code \fhbvmnew\, can be downloaded at the same URL \cite{fhbvm}, where the previous codes are available as well.

The calling sequence of the code is

\begin{equation}\label{call}
{\tt [t,y,etim] = \fhbvmnew(fun,y0,T,M,mu,rho)}
\end{equation}

\medskip
\no In input, with reference to (\ref{fde1}) and (\ref{mixed}):\footnote{Hereafter, in order to have the same calling sequence for both {\tt fhbvm2} and \fhbvmnew, for this latter code, with reference to (\ref{call}): {\tt M} is renamed {\tt N}, {\tt rho} is renamed {\tt n}, and {\tt mu} is renamed {\tt nu}.  As matter of fact, when possible, the two codes will be used with the same input arguments. See \cite{BGIV2025_1} for more details.}
\begin{itemize}
\item {\tt T} is the final time $T$;
\item {\tt M} contains the value $M$ in (\ref{mixed}) (this parameter is named {\tt N}, in the actual interface of the code);
\item {\tt mu} contains $\mu$ (this parameter is named {\tt nu} in the actual interface of the code);
\item {\tt rho} contains $\rho$ (this parameter is named {\tt n} in the actual interface of the code);
\item {\tt y0} is a $\ell \times m$ matrix with the initial conditions (row vectors);
\item {\tt fun} is the identifier of a function implementing the problem (\ref{fde1}). In more detail:
\begin{itemize}
\item {\tt [alpha,m1] = fun()} ~returns $\alpha_1$, $\alpha_2$ and $m_1$. In case $\alpha_1=\alpha_2$ or $\alpha_2=0$, $m_1$ is not used;
\item {\tt fun(t,y)} ~returns the vector field (also in vector mode);
\item {\tt fun(t,y,1)} ~returns the Jacobian matrix of $f$.
\end{itemize}
\end{itemize} 
In output:
\begin{itemize}
\item {\tt t} is the used mesh (\ref{hn});
\item {\tt y} contains the computed approximations (row-wise);
\item {\tt etim} is the elapsed time (in sec).
\end{itemize}

\section{Numerical tests}\label{numtest}

Here, we report some comparisons among the following Matlab\cpr codes, publicly available on the internet:

\begin{itemize}
\item {\tt fde12} \cite{Garr18}, implementing a second-order predictor-corrector method, based on fractional versions of the Adams methods, which can handle only the case $\nu=1$ in (\ref{fde1}). We shall consider the PECE implementation of the methods ({\tt fde12}) and the one with 10 corrector iterations ({\tt fde12-10});
\item {\tt flmm2} \cite{Garr18}, implementing fractional versions of the trapezoidal rule ({\tt flmm2-1}), of the Newton-Gregory formula ({\tt flmm2-2}), and of the BFD2 method ({\tt flmm2-3}). This code can handle only the case $\nu=1$ in (\ref{fde1}), and will be used with the following parameters:
\begin{itemize}
\item tolerance for the nonlinear iteration: {\tt tol=1e-15},
\item maximum number of nonlinear iterations per step: {\tt maxit=1000};
\end{itemize}
\item \FDEpc\, \cite{Garr18}, implementing a second-order product-integration with predictor-corrector, which can handle multi-order FDEs. As with {\tt fde12}, we shall consider the PECE implementation, and that with 10 iterations of the corrector (\FDEpcc);
\item \FDEim\, \cite{Garr18}, implementing a second-order implicit product-integration of trapezoidal type, which can handle multi-order FDEs. It is used with the same parameters {\tt tol} and {\tt itmax} considered for {\tt flmm2};
\item {\tt fhbvm} \cite{BGI2024}, based on a FHBVM(22,20) method, which can handle only the case $\nu=1$ in (\ref{fde1});
\item  {\tt fhbvm2} \cite{BGIV2025}, based on a FHBVM(22,22) method, which can handle only the case $\nu=1$ in (\ref{fde1});
\item the code \fhbvmnew\,here described, based on {\blue either a FHBVM(22,22) or a FHBVM(30,22) method, as explained above}, which can handle the cases $\nu=1$ and $\nu=2$ in (\ref{fde1}). In the case $\nu=1$, it will be used with the same parameters used for {\tt fhbvm2}.
\end{itemize}
We refer to \cite{BGIV2025_1} for a quick introduction to the codes {\tt fde12}, {\tt flmm2}, {\tt fhbvm}, and  {\tt fhbvm2}. Moreover, the code \FDEpc\, is 
a kind of multi-order extension of {\tt fde12}, and \FDEim\, is a kind of multi-order extension of {\tt flmm2}, both couples of codes having a similar interface.\footnote{These codes are available at the URL \cite{garso}.}

All numerical tests have been done on a  12-core M4-pro Silicon based computer with 64GB of memory, using Matlab\cpr Rel.\,2025a.

When necessary, comparisons among the codes, on a given test problem, will be done by using a corresponding {\em work-precision diagram (WPD)}, where the execution time (measured in sec) is plotted againts accuracy, measured in terms of {\em mixed-error significant computed digits (mescd)}, defined as follows:
\begin{equation}\label{mescd}\blue
\mbox{mescd}~=~\max\left\{0,\,-\log_{10}\max_n \|(y(t_n)-y_n)./(1+|y(t_n)|)\|\right\},
\end{equation}  
with $y_n\approx y(t_n)$, this latter the reference solution, $|\cdot|$ the vector of the absolute values, and $./$ the component-wise division.
As is clear, higher and higher accuracies are obtained by suitably tuning the parameters of the various codes. 

The problems with $\nu=1$ are taken from the FDE-Testset \cite{BGIV2025_1}: in particular, at the URL \cite{testset} Matlab\cpr codes implementing them, and providing reference solutions, are available.

{\blue
\begin{rem} For the sake of completeness, we recall that in the definition (\ref{mescd}) the absolute error
$|y(t_n)-y_n|$ is scaled, component-wise, by $(1+|y(t_n)|)$, so that each entry of the vector
 $|y(t_n)-y_n|./(1+|y(t_n)|)$ is approximately the absolute error, if the corresponding entry of $|y(t_n)|$ is much smaller than 1, or the relative error,  if much larger that 1. This measure is widely used in the numerical comparison of numerical codes (see, e.g., the Test-Set for IVP Solvers \cite{ivpsolv}). 
\end{rem}
}

\begin{figure}[t]
\includegraphics[width=16cm]{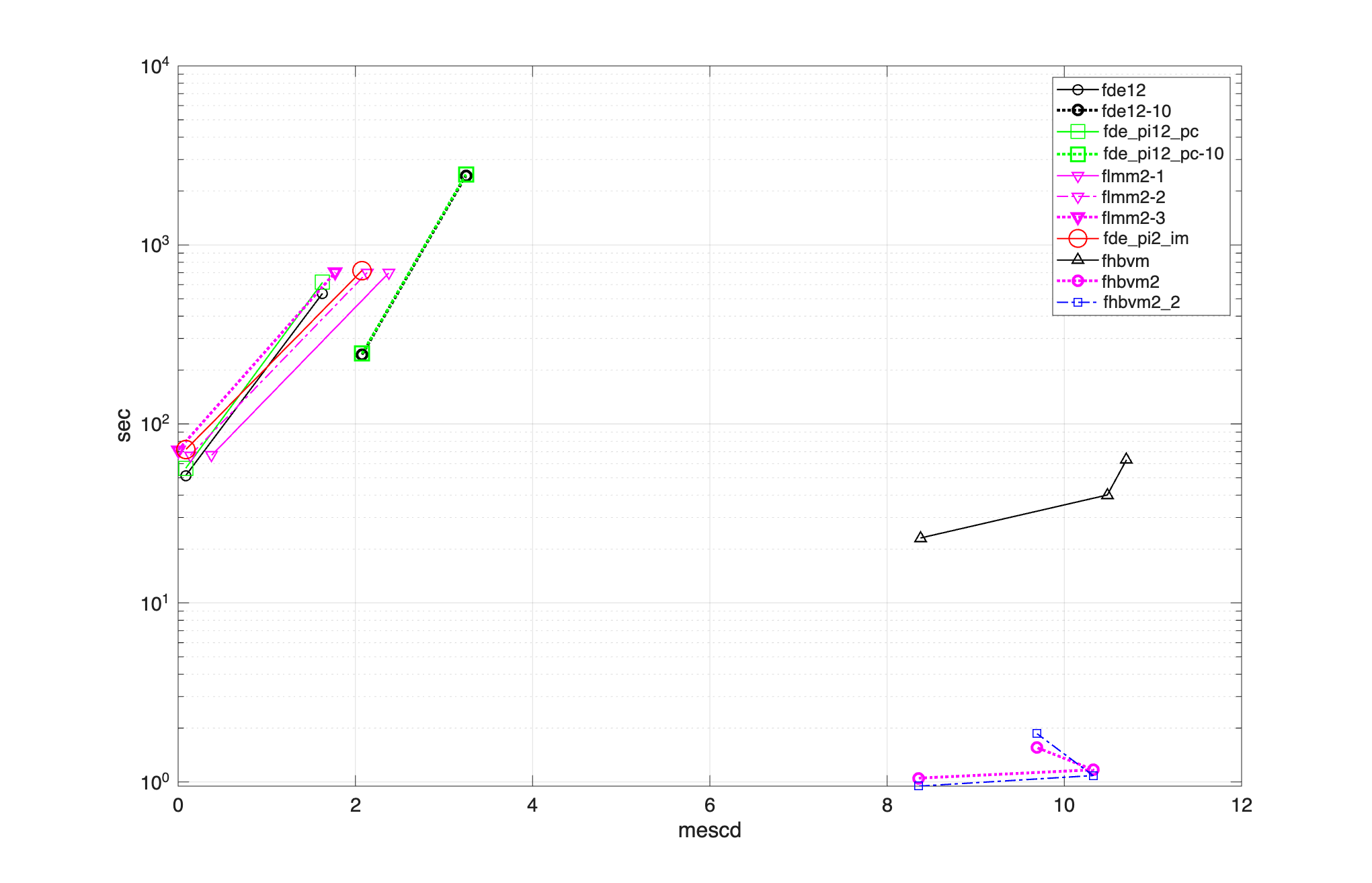} 
\caption{WPD for Problem (\ref{prob2}).}
\label{prob2fig}
\end{figure}

\subsection*{\blue Problem~1} 
This is a {\em stiffly oscillatory} linear problem with only one fractional derivative. It is taken from \cite[Problem~4]{BGIV2025_1} and, at the URL \cite{testset}, a corresponding Matlab\cpr code, implementing the problem and providing the reference solution, is available:
\begin{equation}\label{prob2}
y^{(0.5)} = \frac{1}8\pmatrix{rrrrr}
41 &41 &-38 &40 &-2\\
-79 &81 &2 &0 &-2\\
20 &-60 &20 &-20 &-8\\
-22 &58 &-24 &20 &-4\\
1 &1&-2 &-4 &-2
 \endpmatrix y, \qquad t\in[0,20], \qquad y(0) = \pmatrix{r}1\\ 2\\ 3\\ 4\\ 5\endpmatrix.\end{equation}
Based on \cite{BGIV2025_1}, the following parameters are used for the various codes, to construct the corresponding WPD:
\begin{itemize}
\item {\tt fde12}, {\tt fde12-10}, \FDEpc, \FDEpcc: ${\tt h\,} = 10^{-i}$, $i=5,6$;
\item  {\tt flmm2-1}, {\tt flmm2-2}, {\tt flmm2-3}, \FDEim: ${\tt h\,} = 10^{-i}$, $i=4,5$;
\item {\tt fhbvm}: ${\tt M\,} = 100i$, $i=3,4,5$;
\item {\tt fhbvm2}, \fhbvmnew: ${\tt nu\,}=50$, ${\tt n\,} = 1$, ${\tt N\,} = 100i$, $i=3,4,5$.
\end{itemize}
The corresponding WPD is reported in Figure~\ref{prob2fig}, from which one obtains that:
\begin{itemize}
\item {\tt fde12} and \FDEpc\, reach an accuracy of less than 2 mescd in more than  500 sec. Increasing the number of corrector iterations improves accuracy to about 3 mescd, but with a much higher execution time (about 2500 sec) for both codes;
\item {\tt flmm2} (all versions) and \FDEim\, can achieve about 2 mescd in about 800 sec;
\item {\tt fhbvm} is the most accurate code, with about 11 mescd, with an execution time of about 60 sec;
\item {\tt fhbvm2} and \fhbvmnew\, do have a comparable performance (over 10 mescd in about 1 sec), as expected, due to the oscillatory nature of the solution \cite{BGIV2025_1}.
\end{itemize}

\subsection*{\blue Problem~2} 
\begin{figure}
\vspace{-1cm}
\centering
\includegraphics[width=16cm]{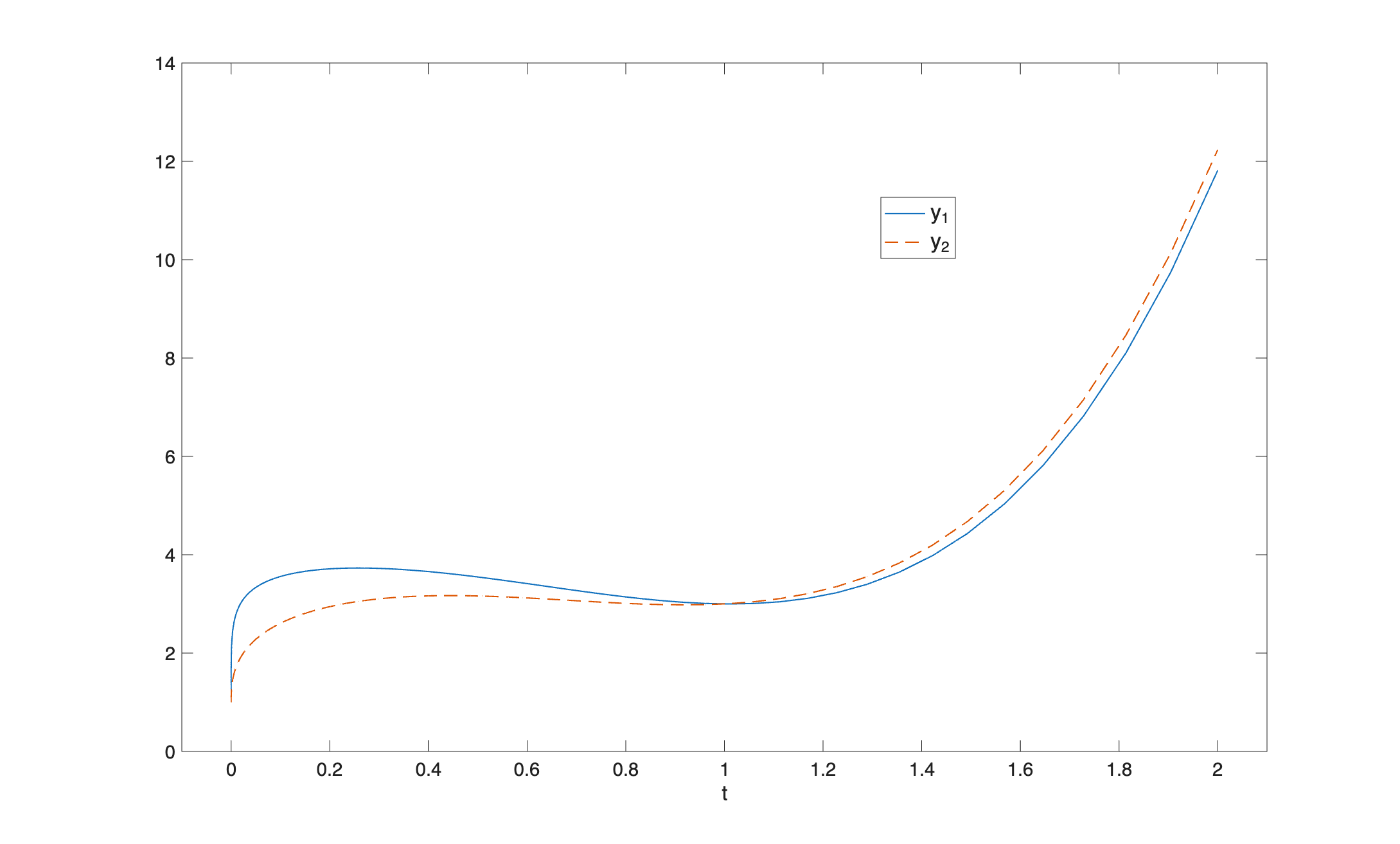}
\includegraphics[width=16cm]{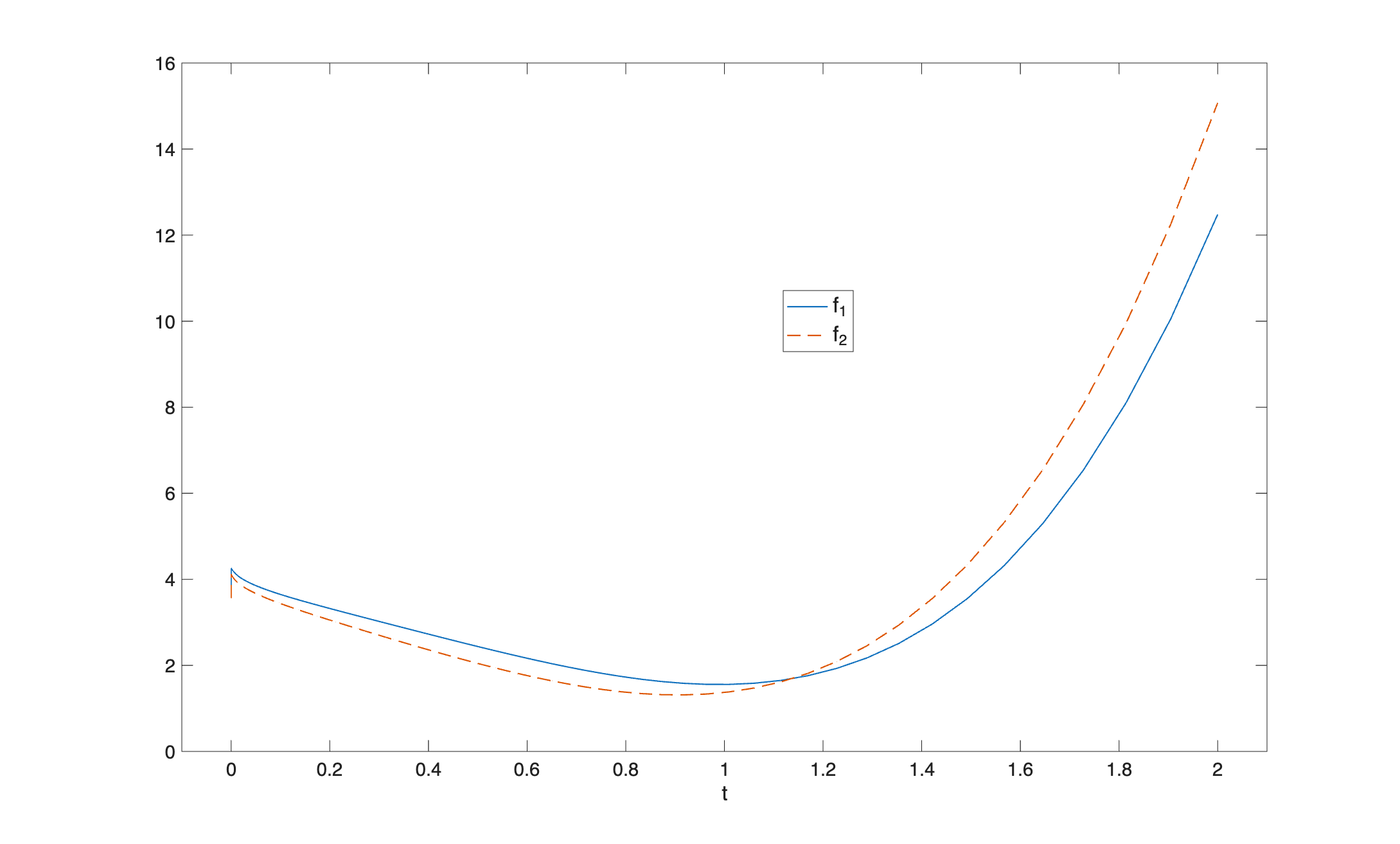}
\caption{Solution (upper-plot) and vector field (lower-plot) for Problem (\ref{prob3}).}
\label{prob3sol}
\end{figure}

The next problem is a multi-order variant of \cite[Problem~7]{BGIV2025_1}:
\begin{eqnarray}\nonumber
y_1^{(\aa_1)} &=& s(t,\aa_2,\beta)^2-y_2^2+g(t,\aa_1,\beta),\\  \nonumber
y_2^{(\aa_2)} &=&  -s(t,\aa_1,\beta)^2+y_1^2+g(t,\aa_2,\beta),\qquad t\in[0,2],\\ \label{prob3}
y_1(0) &=& y_2(0)~=~1, 
\end{eqnarray}
having solution ~$y_i(t) = s(t,\aa_i,\beta)$, $i=1,2$, with
$$s(t,\aa,\beta) = (1-t^2)^2 +4t^\aa +(2-3t^{0.2})t^{\aa+\beta},$$ 
and
$$g(t,\aa,\beta) = \frac{24 t^{4-\aa}}{\Gamma(5-\aa)} - \frac{4 t^{2-\aa}}{\Gamma(3-\aa)}
-3t^{0.2+\beta}\frac{\Gamma(1.2+\aa+\beta)}{\Gamma(1.2+\beta)}
+2t^\beta\frac{\Gamma(1+\aa+\beta)}{\Gamma(1+\beta)}+4\Gamma(1+\aa).
$$
We shall consider the values ~$\aa_1=0.2$, ~$\aa_2=0.4$, ~$\beta=0.1$. The solution, and the corresponding vector field, are depicted in Figure~\ref{prob3sol}.

\begin{figure}
\vspace{-2.5cm}
\centering     
 \includegraphics[width=16cm]{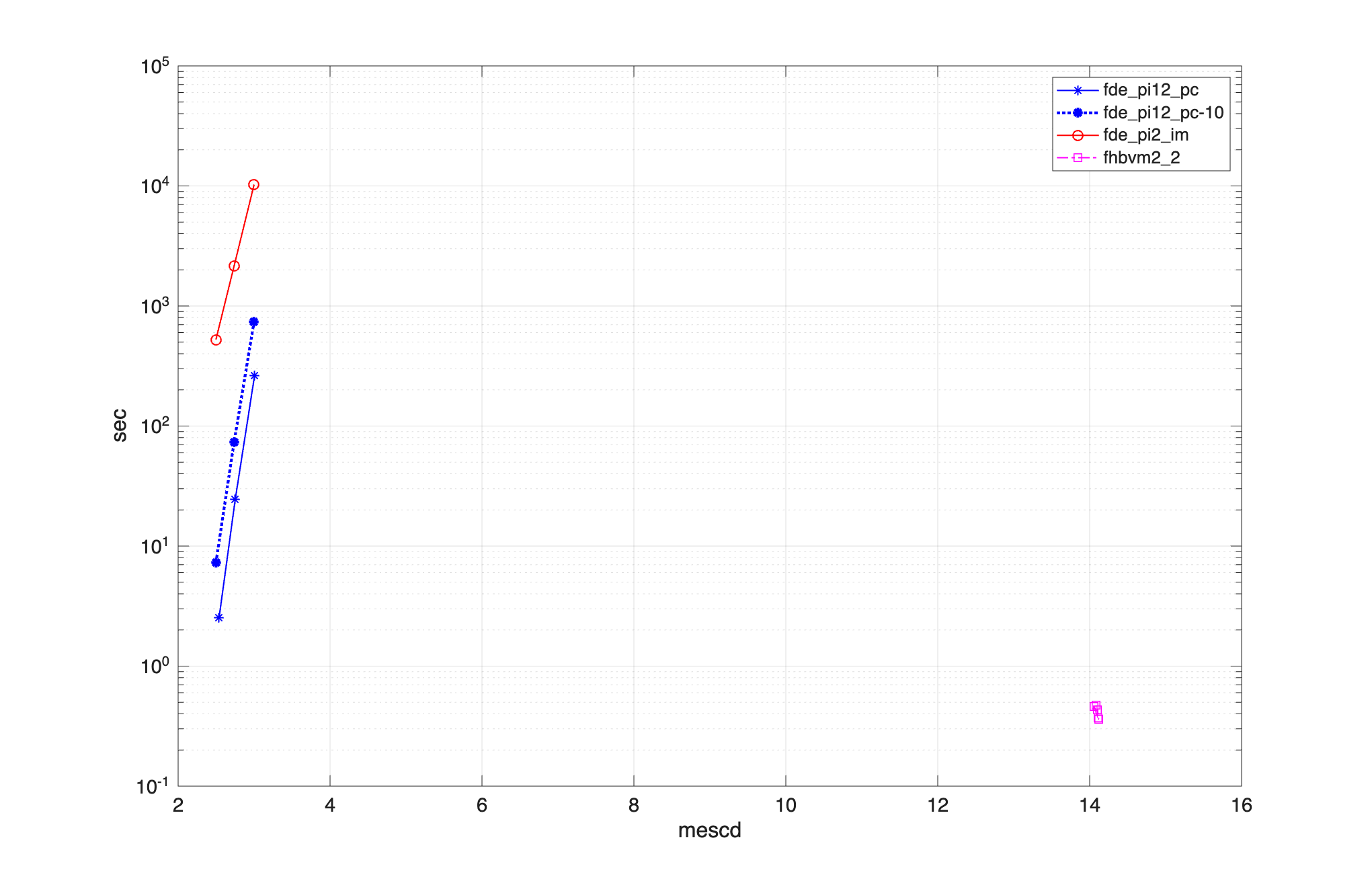} 
\caption{WPD for Problem (\ref{prob3}).}
\label{prob3fig}
 \includegraphics[width=16cm]{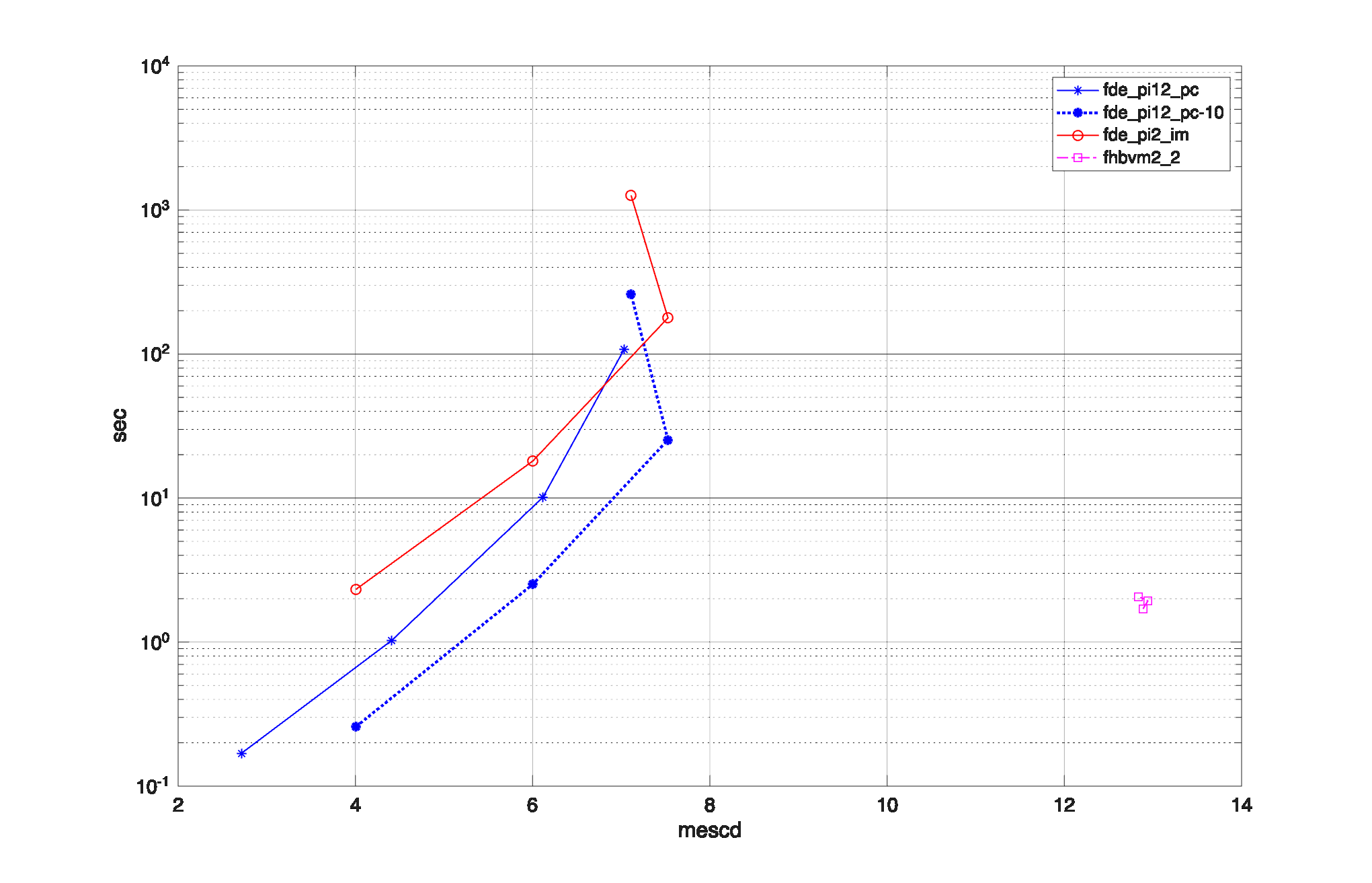}
\caption{WPD for Problem (\ref{prob4}), $\aa_1=0.7$, $\aa_2=0.8$.}
\label{prob4fig}
\end{figure}

The following parameters are used for the codes {\tt fde12}, {\tt flmm2}, and \fhbvmnew, which are the only ones that can be used, in this case:
\begin{itemize}
\item \FDEpc, \FDEpcc,  \FDEim:  ${\tt h\,} = 10^{-i}$, $i=5,6,7$;
\item \fhbvmnew: ${\tt n\,} = 2$, ${\tt nu\,} = 100$, ${\tt N\,} = 5i$, $i=2,3,4,5,6$.
\end{itemize}
The obtained results are reported in the WPD in Figure~\ref{prob3fig}, from which one may infer that:
\begin{itemize}
\item \FDEpc\, and \FDEpcc\, can reach an accuracy of about 3 mescd in about 250 sec and 700 sec, respectively. Consequently, increasing the corrector iterations does not improve accuracy but only affects the execution time;
\item \FDEim\, reaches the same accuracy of about 3 mescd in about $10^4$ sec;
\item \fhbvmnew\, reaches full machine accuracy (more than 14 mescd) in about 0.3 sec.
\end{itemize}

{\blue
We use problem (\ref{prob3}) also for assessing the convergence order (\ref{finerr}), when using a graded mesh (\ref{graded}). For this purpose, we use  a modification of the code \fhbvmnew\, implementing a FHBVM$(30,s)$ method. We use the parameters ${\tt n} = {\tt N} = N$, so that $h_N=2/N$, whereas the parameter {\tt nu} is chosen so that $h_N^{s+\aa}=h_1^{2\aa}$. In fact, in such case (see (\ref{talfa})) $\taa = \aa=0.2$, and the upper bound in (\ref{finerr}) (graded case) turns out to be optimal. In doing this, we also consider that the parameter $r$ defining the graded mesh is given by (\ref{rmixed}), i.e., $r=N/(N-1)$. As one may see, the numerical results listed in Table~\ref{prob3tab} are in good agreement with (\ref{finerr}), since the convergence order appears to be approximately given by $s+0.2$, as expected, for the considered values $s=2,3,4$.
\begin{table}[t]\blue
\caption{Convergence results for HBVM$(30,s)$ used with a graded mesh with $h_N = 2/N$.}
\label{prob3tab}
\centering
\begin{tabular}{|r|rr|rr|rr|}
\hline
\hline
 & \multicolumn{2}{c|}{$s=2$}   &  \multicolumn{2}{c|}{$s=3$} & \multicolumn{2}{c|}{$s=4$}\\
 \hline
$N$ & error & order & error & order & error & order\\
\hline  \hline
10 & 2.88e-02 &   --- & 4.36e-04 &   --- & 2.07e-06 &   --- \\ 
 20 & 6.89e-03 &   2.06 & 5.08e-05 &   3.10 & 1.17e-07 &   4.14 \\ 
 40 & 1.57e-03 &   2.13 & 5.68e-06 &   3.16 & 6.41e-09 &   4.19 \\ 
 80 & 3.48e-04 &   2.17 & 6.21e-07 &   3.19 & 3.43e-10 &   4.22 \\ 
 \hline
 \hline
\end{tabular}
\end{table}

}

\subsection*{\blue Problem~3} 
\begin{figure}
\centering
\vspace{-1cm}
\includegraphics[width=16cm]{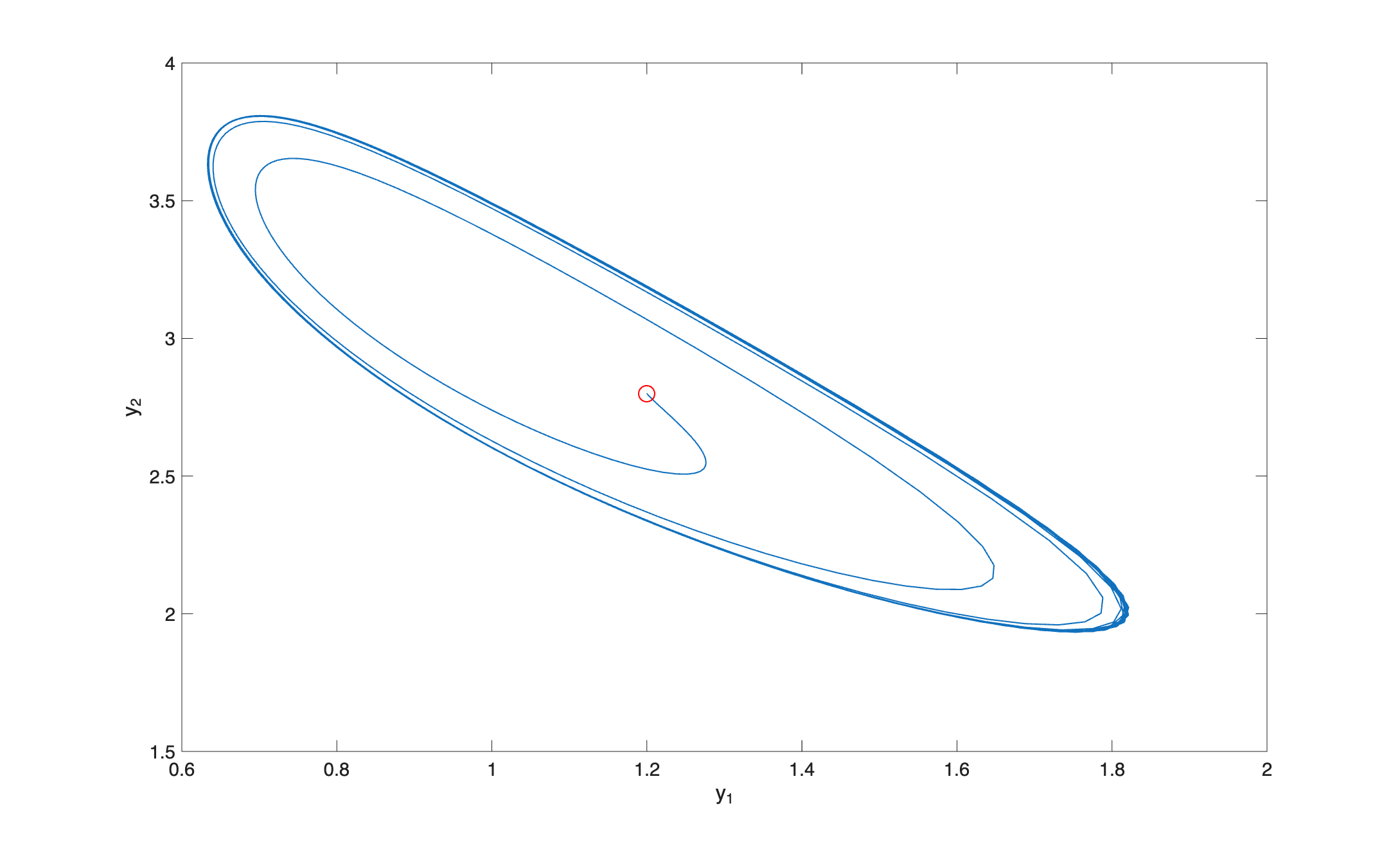}
\includegraphics[width=16cm]{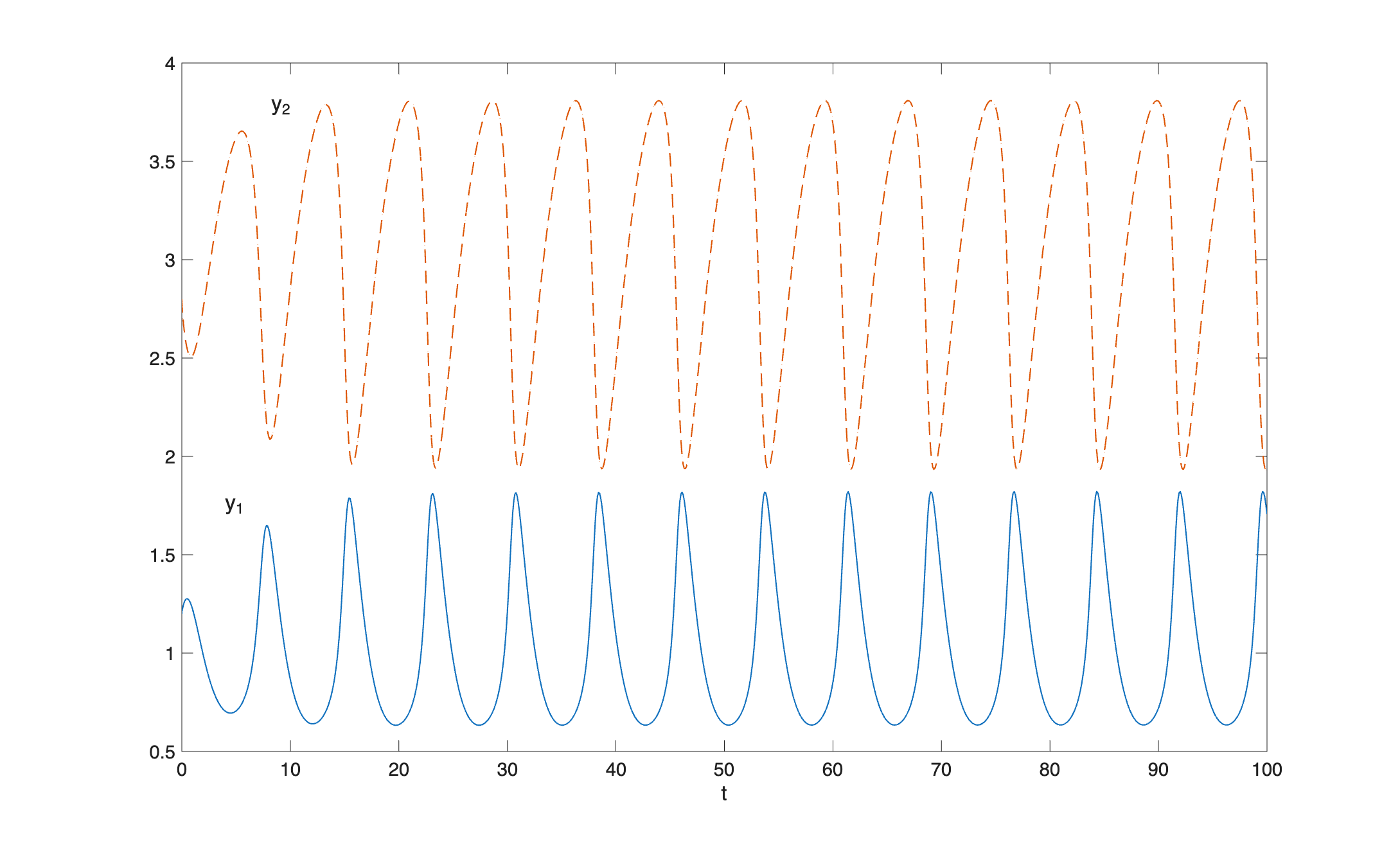}
\caption{Solution in the phase space (upper-plot) and versus time (lower-plot) for Problem (\ref{prob4}).}
\label{prob4sol}
\end{figure}

This problem is a multi-order fractional version of the Brusselator problem, adapted from \cite[Problem~(32)]{Garr18}:
\begin{eqnarray}\nonumber
y_1^{(\aa_1)} &=& A-(B+1)y_1+y_1^2y_2,\\ \nonumber
y_2^{(\aa_2)} &=&  By_1-y_1^2y_2,\qquad\qquad\qquad t\in[0,100],\\ \label{prob4}
y_1(0) &=& 1.2, \qquad y_2(0)~=~2.8, 
\end{eqnarray}
where we use the same parameters as in \cite{Garr18}: $A=1$, $B=3$, $\aa_1=0.8$, and $\aa_2=0.7$.
The solution is depicted in Figure~\ref{prob4sol}, both in the phase space and versus time, from which one infers that it approaches a limit cicle, thus becoming of periodic type: in the first plot, the circle denotes the initial point of the trajectory.

\no The following parameters are used for the codes \FDEpc, \FDEim, and \fhbvmnew, which are the only ones that can be used, in this case:
\begin{itemize}
\item \FDEpc, \FDEpcc, \FDEim: ${\tt h\,} = 10^{-i}$, $i=2,3,4,5$;
\item \fhbvmnew: ${\tt n\,} = 1$, ${\tt nu\,} = 50$, ${\tt N\,} = 50i$, $i=4,5,6$.
\end{itemize}
The obtained results are reported in Figure~\ref{prob4fig}, where the accuracy w.r.t. a reference solution at the end point is reported. This latter has been evaluated numerically by using \fhbvmnew\, on 3 consecutive doubled meshes, resulting into:  
$$y_1(100) \approx 1.706502172199, \qquad y_2(100) \approx 1.940414058005.$$
From the obtained results, one concludes that:
\begin{itemize}
\item \FDEpcc\, reaches an accuracy of 7.5 mescd in about 25 sec, whereas \FDEpc\, reaches about 7 mescd in 100 sec; 
\item \FDEim\, reaches an accuracy of 7.5 sec in approximately 180 sec;
\item \fhbvmnew\, reaches 13 mescd in less than 2 sec, independently of the stepsize used, as one expects from a spectrally accurate method.
\end{itemize}

\subsection*{\blue Problem~4}  
\begin{figure}[t]
\includegraphics[width=16cm]{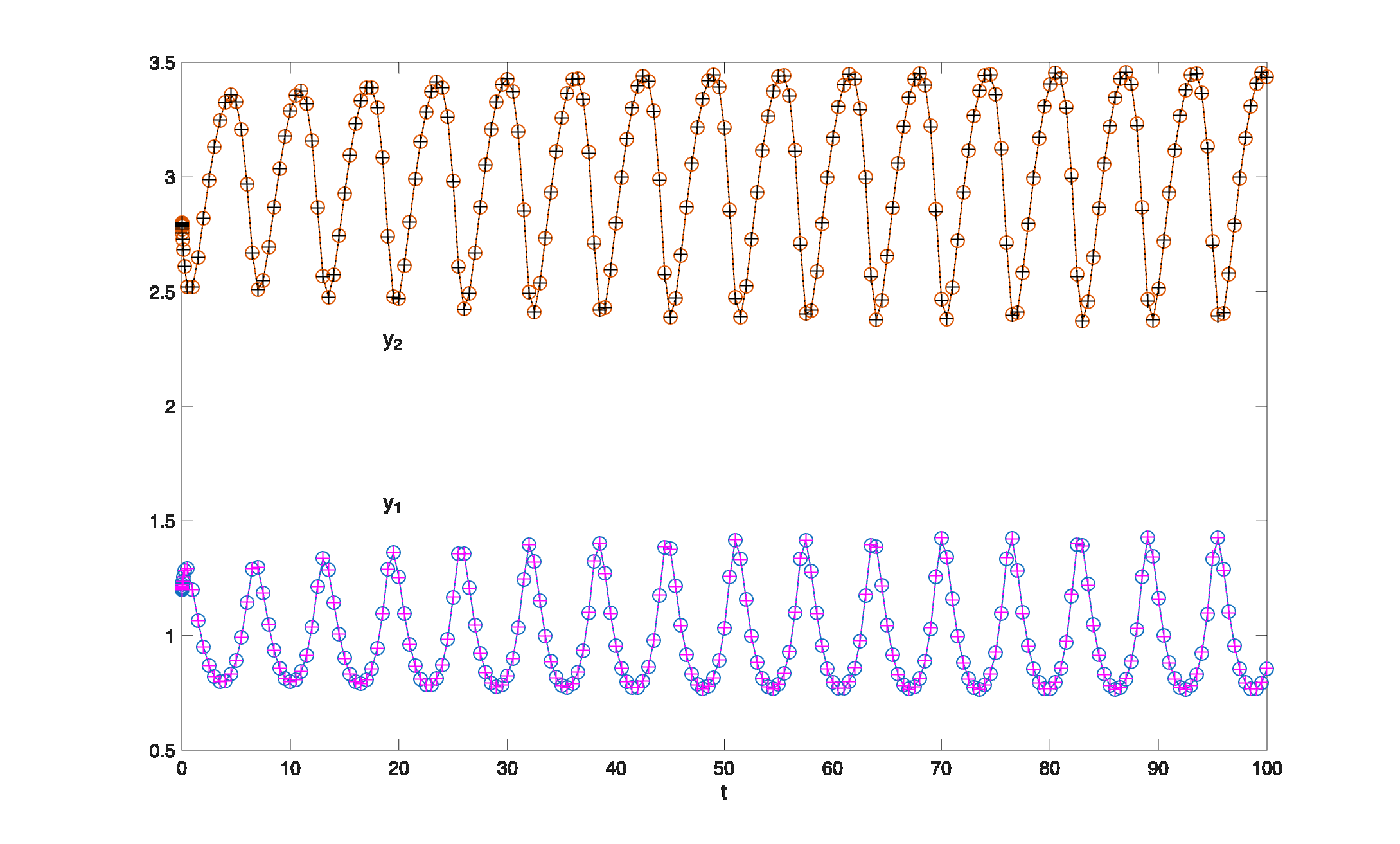}
\caption{Solution of Problem (\ref{prob4}), $\aa_1=\aa_2=0.7$ (circles) and $\aa_1=0.7$, $\aa_2=0.7+10^{-4}$ (crosses).}
\label{prob5fig}
\end{figure}

As stressed in Section~\ref{nu_1}, the {\em blended iteration}, which can be used in the case $\nu=1$ in (\ref{fde1}), only requires to factor one matrix having the same size as that of the continuous problem. It has, therefore, a lower computational cost, w.r.t. the simplified Newton iteration, to be used in the case $\nu>1$. The next example is aimed at assessing this fact. For this purpose, we shall use the code \fhbvmnew\, for solving problem (\ref{prob4}) with the same parameters $A=1$ and $B=3$  used before, but with $\aa_1$, and $\aa_2$ chosen as follows:
\begin{itemize}
\item $\aa_1=\aa_2=0.7$;
\item $\aa_1=0.7$, ~$\aa_2=0.7+\epsilon$, ~with~ $\epsilon=10^{-4}$.
\end{itemize}
As is clear, in the first case only one fractional derivative order occurs (i.e., $\nu=1$ in (\ref{fde1})), whereas they are two in the second case (i.e., $\nu=2$ in (\ref{fde1})). Despite this fact,  the computed solutions are very similar, as is shown in Figure~\ref{prob5fig}. Consequently, we shall assume that the code provides a comparable accuracy, when used with the same parameters, in both cases. In particular, similarly as done for Problem~4, we shall use the parameters ${\tt n\,} = 1$, and ${\tt nu\,}=50$. Instead, the values of {\tt N} vary as reported in Table~\ref{prob5tab}, where we also list:
\begin{itemize}
\item the required number of fixed-point iterations (\ref{fixit}) (basically, used in the initial graded mesh, where the timesteps are very small);
\item the number of blended iterations (\ref{blend}), in the case $\nu=1$, or simplified Newton iterations (\ref{simpNewt1}), in the case $\nu=2$; 
\item the execution times (in sec).
\end{itemize}
As expected, the execution times, when using the blended iteration, are lower than those required by the simplified Newton iteration, despite the fact that the latter requires less steps to reach convergence (approximately 30\%  less). In particular, when using the blended iteration,  the code is more than 3.5 times faster than when using the simplified Newton iteration.

\begin{table}[t]
\caption{Execution times (in sec) and nonlinear iterations  for \fhbvmnew, when used for solving (\ref{prob4}) with $\aa_1=\aa_2=0.7$ (case~1) and $\aa_1=0.7$, $\aa_2=0.7+10^{-4}$ (case~2).}
\label{prob5tab}
\centering\medskip

\begin{tabular}{|c|ccc|ccc|}
\hline
\hline
          &\multicolumn{3}{c|}{case~1} & \multicolumn{3}{c|}{case~2}\\
          \hline
 {\tt N}& fixed-point  & blended  & execution & fixed-point  & Newton & execution\\
          & iterations    & iterations & time & iterations   & iterations & time\\
\hline
200 &  150 & 2679 & 0.40 &  142 & 1950 & 1.51 \\ 
250 &  149 & 3193 & 0.45 &  141 & 2238 & 1.65 \\ 
300 &  146 & 3644 &  0.54 & 146 & 2529 &  2.05 \\ 
350 &  146 & 4115 &  0.64 &  146 & 2822 &  2.37 \\ 
400 &  152 & 4627 &  0.73 &  144 & 3084 &  2.80 \\ 
450 &  151 & 5116 &  0.85 &  143 & 3368 &  3.14 \\ 
500 &  151 & 5617 &  0.99 &  143 & 3667 &  3.54 \\ 
550 &  150 & 5969 &  1.07 &  150 & 3934 &  3.83 \\ 
600 &  148 & 6423 &  1.19 &  148 & 4146 &  4.22 \\ 
\hline
\hline
\end{tabular}

\end{table} 

\subsection*{\blue Problem~5} 

\begin{figure}
\vspace{-2cm}
\centering
\includegraphics[width=15cm,height=10cm]{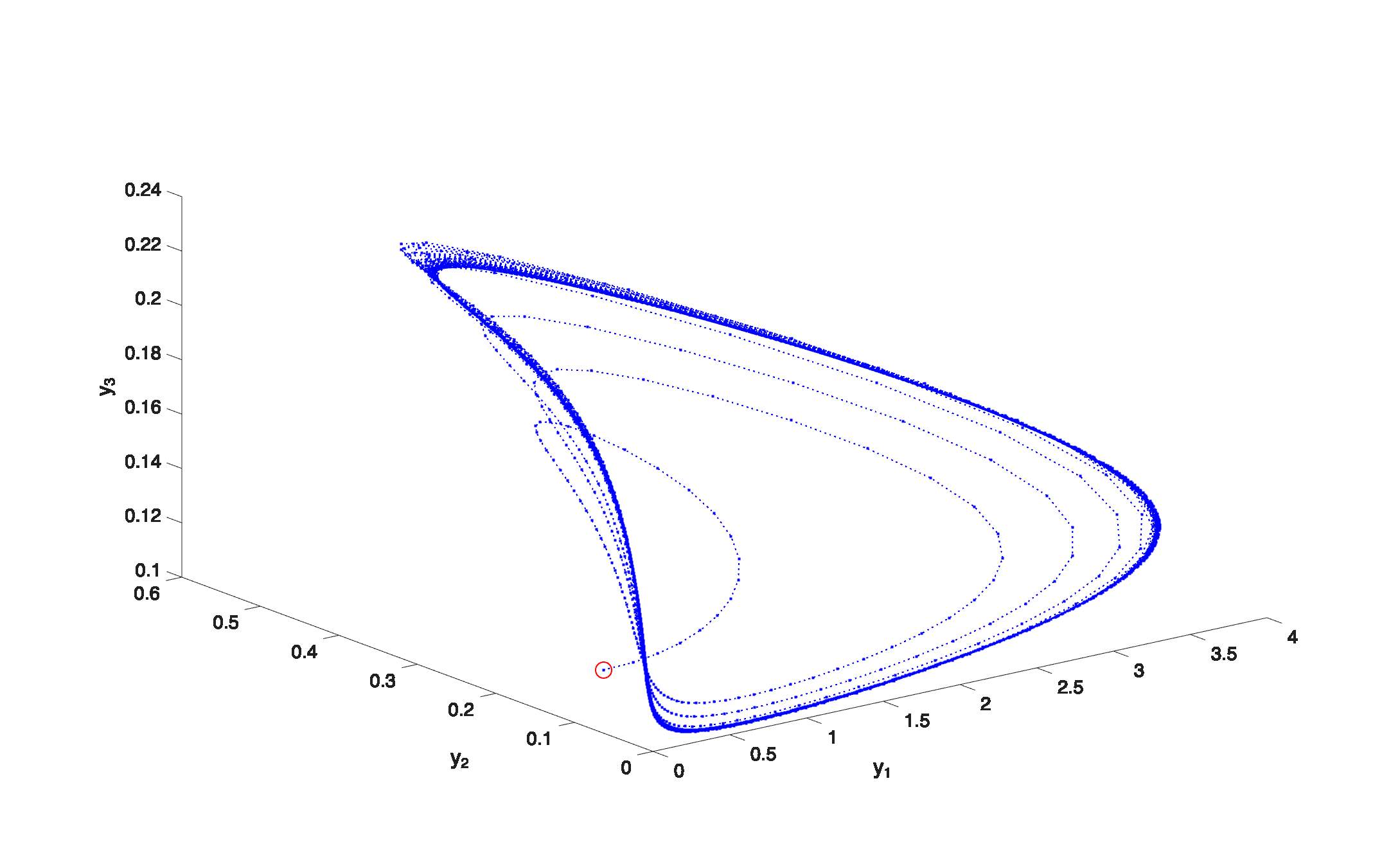}
\includegraphics[width=16cm]{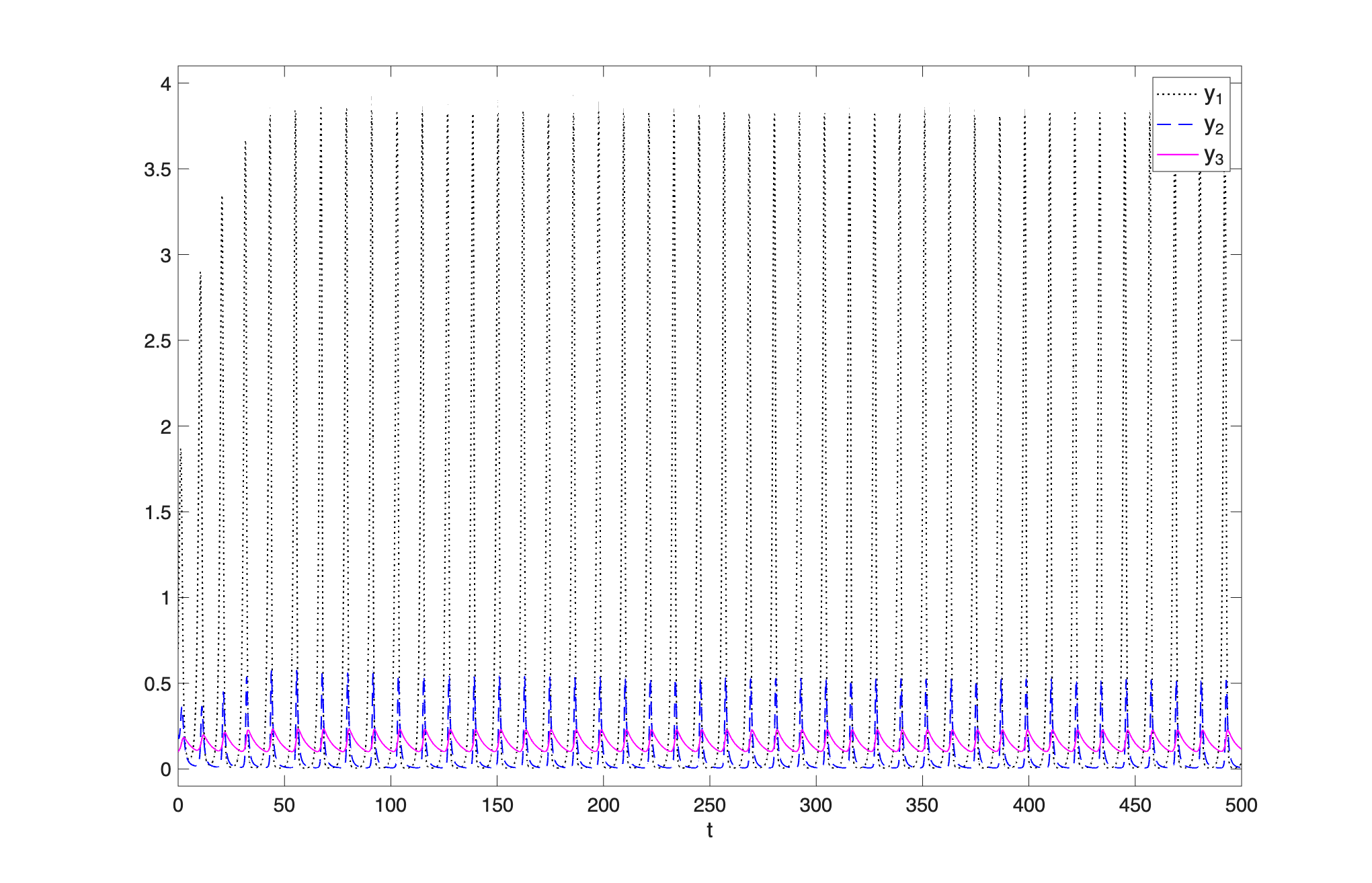}
\caption{Solution of Problem (\ref{prob6}) in the phase space (upper plot) and versus time (lower plot). In the upper plot, the circle denotes the initial point of the trajectory.}
\label{prob6fig}
\end{figure}

At last, we consider the following multi-order fractional version of the predator-prey model described in \cite{Panja2019}, including intra-species competition:
\begin{eqnarray}\nonumber
y_1^{(\aa_1)} &=& r_1y_1-a_{11}y_1^2-a_{12}y_1y_2-a_{13}y_1y_3,\\[2mm] \nonumber
y_2^{(\aa_2)} &=& a_{21}y_1y_2 -a_{22}y_2^2 -\frac{a_{23}y_2y_3}{1+\beta y_2}-r_2y_2, \qquad\qquad t\in[0,T],\\[2mm] \label{prob6}
y_3^{(\aa_2)} &=& a_{31}y_1y_3+\frac{a_{32}y_2y_3}{1+\beta y_2}-a_{33}y_3^2-r_3y_3, \qquad\qquad y(0) ~=~ y^0\in\RR^3.
\end{eqnarray}
Here:
\begin{itemize}
\item $y_1$ represents the density of the preys;
\item $y_2$ is the density of the {\em intermediate} predators, which cannot prey the {\em top} predators;
\item $y_3$ denotes the density of the {\em top} predators, which can prey also the {\em intermediate} ones.
\end{itemize}
We use the following parameters: ~$\aa_1 = 0.99$, $\aa_2=0.8$; $r_1=5$, $r_2=1$, $r_3=0.1$; 
$a_{11} = 0.01$, $a_{12}=1$, $a_{13}=35$,  $a_{21}=1$, $a_{22}=0.2$, $a_{23}=1$,  $a_{31}=0.1$,  $a_{32}=1$,  $a_{33}=0.3$;
$\beta=0.01$. Moreover, we use the initial point $y^0=\left(0.7, \, 0.2,\, 0.1\right)^\top$, and integrate up to $T=500$: as one may infer from the plots in Figure~\ref{prob6fig}, the solution tends to become of periodic type, with a period $\tau\approx11.8$. Since a reference solution is not available for this problem, to obtain corresponding accuracy results we run each code on consecutive doubled meshes, as below specified:
\begin{itemize}
\item \FDEpc, \FDEpcc, \FDEim: ${\tt h\,} = 10^{-2}\cdot 2^{1-\ell}$, $\ell=1,\dots,7$;
\item \fhbvmnew: ${\tt n\,} = 1$, ${\tt nu\,} = 50$, ${\tt N\,} = 500\cdot 2^{\ell-1}$, $\ell=1,2,3,4$.\footnote{Consequently, the timestep in the uniform part of the mesh is\, $h=2^{1-\ell}$, $\ell=1,2,3,4$.}
\end{itemize}
Then, for each code, the accuracy on the $\ell$th mesh is estimated through the solution on the subsequent $(\ell+1)$st mesh.
The obtained results are listed in Table~\ref{prob6tab}, from which one may conclude that:
\begin{itemize}
\item \FDEpc\, reaches 3.5 mescd in about 24 sec;
\item \FDEpcc\,  and \FDEim\, do have a comparable accuracy, and can reach a higher accuracy (slightly larger than 4 mescd) but with a higher execution time. Moreover, that of \FDEim\, is much larger than that of \FDEpcc; 
\item \fhbvmnew\, reaches a very high accuracy (almost 12 mescd) in about 22 sec. Moreover, it is able to achieve more than 10 mescd by using a timestep as large as 1, in the uniform part of the mesh: this reflects the capability, of the underlying FHBVM(22,22) method, of obtaining a spectrally accurate solution. 
\end{itemize}

\begin{table}
\caption{Estimated accuracy (in mescd) and execution times when solving (\ref{prob6}) on consecutive doubled meshes (see text for details).}
\label{prob6tab}
\medskip
\centering

\begin{tabular}{|c|rr|rr|rr|rr|}  
\hline
\hline
          &\multicolumn{2}{c|}{\FDEpc}&\multicolumn{2}{c|}{\FDEpcc}&\multicolumn{2}{c|}{\FDEim}&\multicolumn{2}{c|}{\fhbvmnew}\\
$\ell$ & mescd & sec & mescd & sec & mescd & sec & mescd & sec\\
\hline
 1 &   0.59 &   0.76  &   1.83 &   2.21  &   1.83 &  21.41  &  10.22 &   4.41   \\ 
  2 &   1.16 &   1.54  &   2.43 &   4.40  &   2.43 &  40.73  &  11.35 &   9.29   \\ 
  3 &   1.77 &   3.13  &   3.03 &   8.91  &   3.03 &  78.52  &  11.68 &  22.50   \\ 
  4 &   2.40 &   6.30  &   3.64 &  18.15  &   3.64 & 161.37  &     &       \\ 
  5 &   3.09 &  12.12  &   4.22 &  36.71  &   4.22 & 309.53  &     &       \\ 
  6 &   3.55 &  24.16  &   4.06 &  74.04  &   4.06 & 596.44  &     &       \\ 
\hline
\hline
\end{tabular}

\end{table} 
\begin{figure}[t]
\includegraphics[width=16cm]{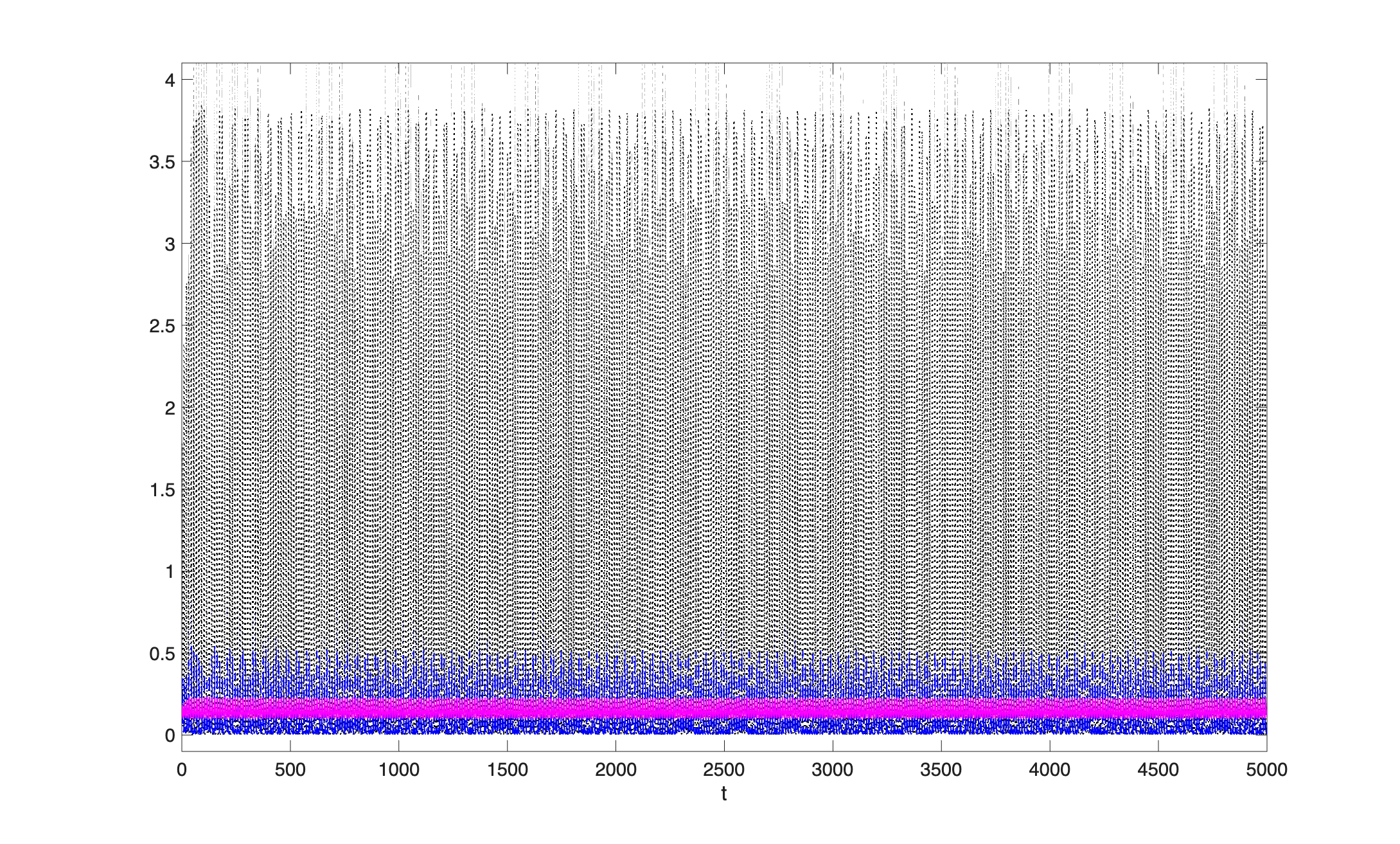}
\caption{\blue Numerical solution of Problem (\ref{prob6}), $T=5000$ by using timestep $h=1$ (except the first 50 initial steps).}
\label{prob6fig1}
\end{figure}

{\blue
\begin{rem}
It is worth mentioning that the considered large width of the integration interval $(T=500)$ of problem (\ref{prob6}) confirms the good stability properties of the FHBVM(30,22) method implemented in the code \fhbvmnew. However, to further confirm this fact, we consider $T=5000$ and, by using the coarsest mesh, among those used for the numerical tests $({\tt n\,} = 1$, ${\tt nu\,} = 50$, ${\tt N\,} = T)$, we obtain the plot in Figure~\ref{prob6fig1}, which still exhibits the correct periodic behavior.
\end{rem}
}

\section{Conclusions}\label{fine}

In this paper, we have extended the class of methods named {\em Fractional HBVMs (FHBVMs)} for numerically solving multi-order FDEs. The implementation details of the methods have been given in the general case: in particular, a comprehensive review of Jacobi-Pin\~eiro quadratures has been presented, with the aim of providing algorithmic insights. As a result, a new Matlab\cpr code, named \fhbvmnew, has been released. The code implements the case of 2 possible different fractional orders, and relies on FHBVM(22,22): this latter method had already been used in the code {\tt fhbvm2}, handling the single-order case. Alike the latter code, also the new one  is able to gain spectrally accurate solutions, on a suitably chosen mesh. Numerical tests taken from the literature duly confirm this fact, also showing its superior performance over existing numerical codes, designed for solving multi-order FDEs, available on the internet. {\blue We plan to extend the code \fhbvmnew\, to the case of more fractional orders, as soon as an effective code for computing the Jacobi-Pi\~{n}eiro weights and abscissae, in the general case, will be available.}

\paragraph*{Data availability.} The Matlab\cpr code \fhbvmnew\, is available at the URL \cite{fhbvm}.
 
\paragraph*{\bf Acknowledgements.} The first three authors are members of the ``Gruppo Nazionale per il Calcolo Scientifico-Istituto Nazionale di Alta Matematica (GNCS-INdAM)''. The last author is supported by the project n.ro PUTJD1275 of the Estonian Research Council. The paper has been written during a period of visit of the last author at the University of Florence.


\paragraph*{Declarations.} The authors declare no conflict of interests.

\end{document}